\documentclass[12pt]{article}
\usepackage{latexsym}
\usepackage{epsfig,graphics}
\def\Re{{\bf R}}

\newtheorem{theorem}{Theorem}[section]

\newtheorem{defin}{Definition}[section]

\begin{document}
\baselineskip=14pt
\title{\large \bf Deterministic global optimization using space-filling curves and
 multiple estimates of   Lipschitz and H\"{o}lder constants}
\author{\small Daniela Lera\thanks{Dipartimento di Matematica e Informatica,
Universit\`a di Cagliari, Cagliari,
Italy}\,\,\, and Yaroslav D. Sergeyev\thanks{Corresponding author;
Dipartimento di Ingegneria Informatica, Modellistica, Elettronica e
Sistemistica, Universit\`{a} della Calabria  and the Institute of
High Performance
  Computing and Networking of the National Research Council of Italy, Via Pietro Bucci
42C, 87036 Rende (CS), Italy; and Software Department, N.I.
Lobachevskiy University of Nizhni Novgorod, Gagarin Av. 23, Nizhni
Novgorod, Russia}}
\date{}
\maketitle

\begin{abstract}
In this paper,    the global optimization problem $\min_{y\in S}
F(y)$ with $S$ being a hyperinterval in $\Re^N$ and $F(y)$
satisfying the Lipschitz condition with an unknown Lipschitz
constant is considered. It is supposed that the function $F(y)$ can
be multiextremal, non-differentiable, and given as a `black-box'. To
attack the problem, a new global optimization algorithm based on the
following two ideas is proposed and studied both theoretically and
numerically. First, the new algorithm uses numerical approximations
to space-filling curves to reduce the original Lipschitz
multi-dimensional problem to a univariate one satisfying the
H\"{o}lder condition. Second, the algorithm  at each iteration
applies a new geometric technique working with a number of possible
H\"{o}lder constants chosen from a set of values varying from zero
to infinity showing so that ideas introduced in a popular DIRECT
method can be used in the H\"{o}lder global optimization.
Convergence conditions of the resulting deterministic global
optimization method are established. Numerical experiments carried
out on several hundreds of test functions show quite a promising
performance of the new algorithm in comparison with its direct
competitors.

 \vspace{10pt} \noindent {\bf Key
Words}. Global optimization, Lipschitz  functions, space-filling
curves,   H\"{o}lder functions, deterministic numerical algorithms,
DIRECT, classes of test functions.

\end{abstract} \vspace{8pt}

\section{Introduction}
\setcounter{equation}{0} Let us consider the following global
optimization problem
\begin{equation} \label{p}
\min\{F(y): \ y\in S=[a,b]\},
\end{equation}
where $[a,b]$ is a  hyperinterval in $\Re^N$. It is supposed that
the objective function $F(y)$ can be   multiextremal, possibly
non-differentiable and it satisfies the Lipschitz condition
\begin{equation} \label{lip}
|F(y')-F(y'')| \leq L \| y'-y''\|, \hspace{1cm} y', y''\in [a,b],
\end{equation}
with an unknown constant $L$, $0<L<\infty,$ in the Euclidean norm.
This statement can very frequently be met in applications
  where each evaluation of $F(y)$ can
be very expensive from the computational point of view (see, e.g.,
\cite{Casado_Garcia_yaro_SIAM,27aa,horst,automation,
pinter,Sergeyev_Daponte,yaro,Zhigljavsky&Zilinskas(2008)}, etc.).
Due to this reason, in the literature  there exist numerous methods
dedicated to the problem (\ref{p}), (\ref{lip}) (together with
references indicated above we can mention such recent publications
as
\cite{Calvin&Zilinskas,Evtushenko&Posypkin,27aa,Kvasov&Sergeyev2009,
Kvasov&Sergeyev2012,PSKJ_2014,Paulavicius&Zilinskas,Zilinskas(2010),
Zilinskas&Zilinskas(2010)}). It should be also mentioned that
methods for problems where the gradient of the objective function
satisfies the Lipschitz condition were also studied (see, e.g.,
\cite{Evtushenko&Posypkin,Gergel_Sergeyev,Kvasov&Sergeyev2015,Lera&Sergeyev2013,Sergeyev_1998,
Sergeyev_Daponte,Sergeyev_Kvasov_book}, etc.).

One of the main issues regarding the problem (\ref{p}), (\ref{lip})
is related to the treatment  of the Lipschitz constant $L$. In the
literature, there exist several approaches for acquiring  the
Lipschitz information that can be distinguished with respect to the
way the Lipschitz constant is estimated during the process of
optimization. For instance, there exist algorithms (see, e.g.,
\cite{Evtushenko&Posypkin,horst,Horst&Tuy,Paulavicius&Zilinskas,piya2})
that use an a priori given estimate $\tilde{L}$ of $L$ (it should be
mentioned that usually in practice it is difficult to obtain valid
estimates) or an adaptive estimate $\tilde{L}_i$ that is
recalculated   at each iteration $i$ of the method (see, e.g.,
\cite{pinter,Sergeyev_Kvasov_book,Sergeyev_Lera_book,yaro}).

 It is
well known that  estimates of the Lipschitz constant have a
significant influence on the convergence and the speed of the global
optimization algorithms (see
\cite{Sergeyev_Kvasov_book,Sergeyev_Lera_book,yaro}). Namely, tight
estimates can accelerate the search, overestimates can slow it down,
and underestimates can lead to losing the global solution.
Algorithms that use in their work a global estimate $\tilde{L}_i$ or
an a priori given estimate $\tilde{L}$  do not take into account any
local information about the behavior of the objective function over
  small subregions of the domain $S$. It has been shown in
\cite{Kvasov&al,daniya,Martinez_Casado,yaro2,Sergeyev_1998,Sergeyev_Lera_book,yaro}
that a smart usage of local estimates $\tilde{L}_i(D_j)$ of the
local Lipschitz constants $\tilde{L}(D_j)$ over subregions $\ D_j
\subset S$ allows one to significantly accelerate the global search.
Naturally, a balancing between the local and global information must
be performed in an appropriate way in order to avoid   missing   the
global solution. Another  interesting approach that has been
introduced in \cite{jones}  uses at each iteration several estimates
of the Lipschitz constant $L$. This algorithm works by partioning
the hyperinterval $S$ in subintervals (called hereinafter also
\textit{subboxes}) and due to this reason it has been called by its
authors \textit{DIRECT} (DIviding RECTangles).

We can see therefore that in the literature there exist four
practical strategies to obtain a Lipschitz information for its
subsequent usage in numerical methods: a)~to consider the real
constant $L$, if it is given, or to use its overestimate  when it is
possible; b)~to calculate dynamically a global a global (i.e., the
same for the whole search region) adaptive estimate of $L$; c)~to
consider estimates of local Lipschitz constants  related to
subintervals of the search region $S$; d)~to take into consideration
a set of   estimates of $L$  selected among all possible values
varying from zero to infinity. In this paper, our attention will be
dedicated to this fourth alternative.

\begin{figure}[t]
\centerline{\epsfxsize=12.5cm\mbox{\epsfbox{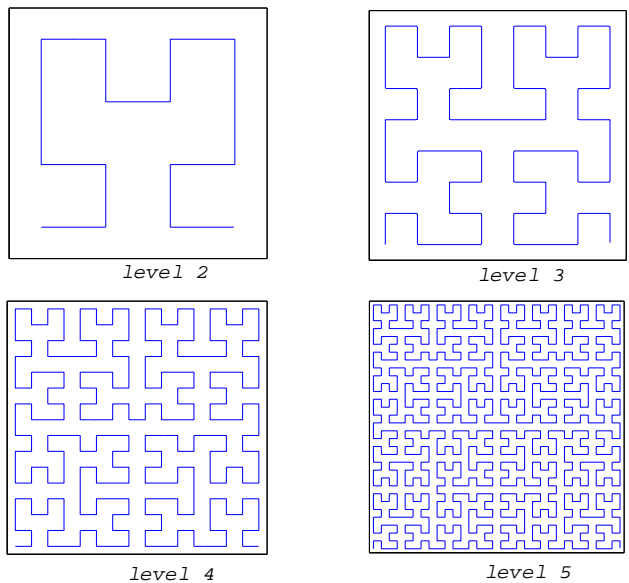}}}
\caption{\em Approximation of the Peano  curve  (Hilbert's version
of its construction is used here) in dimension N=2. } \label{fig.1a}
\end{figure}

As was mentioned above, in the literature there exists a variety of
numerical methods dedicated to the problem (\ref{p}), (\ref{lip}).
One of non trivial views on the problem consists of   the reduction
of the dimension of (\ref{p}), (\ref{lip}) with the help of
space-filling curves. These curves, first introduced by Peano (1890)
and Hilbert (1891) (see
\cite{butz,Sagan,Sergeyev_Lera_book,strongin1,yaro,Strongin_Sergeyev_2003}),
emerge as the limit objects generated by an iterative process. They
are   fractals constructed using the principle of self-similarity.
Peano-Hilbert curves fill in the hypercube $[a,b] \subset \Re^N$,
i.e., they pass through every point of $[a,b]$, and this gave rise
to the term space-filling curves. Examples of construction of these
curves in two dimensions are given in Fig.~\ref{fig.1a}.

It has been shown (see \cite{butz},\cite{strongin1},\cite{yaro})
that, by using space filling curves, the multidimensional global
minimization problem (\ref{p}), (\ref{lip}) can be turned into a
one-dimensional problem. In particular, Strongin has proved  in
\cite{strongin1} that finding the global minimum of the Lipschitz
function $F(y), y \in R^N,$ over a  hypercube is equivalent to
determining the global minimum of the function $f(x)$ in an
interval, i.e., it follows
\begin{equation} \label{fun}
f(x)=F(p(x)),  \hspace{1cm} x\in [0,1],
\end{equation}
where $p(x)$ is the Peano curve. Moreover, the H\"{o}lder condition
\begin{equation} \label{hold}
|f(x')-f(x'')| \leq H |x'-x''|^{1/N}, \hspace{1cm} x',x'' \in
[0,1],
\end{equation}
holds (see \cite{yaro})  for the function $f(x)$ with the constant
\begin{equation} \label{costhold}
H=2L\sqrt{N+3},
\end{equation}
where $L$ is the Lipschitz constant of the original multidimensional
function $F(y)$ from (\ref{p}), (\ref{lip}).

Thus, one can  try to solve the problem (\ref{p}), (\ref{lip}) by
using algorithms proposed for minimizing H\"{o}lderian functions
(\ref{fun}), (\ref{hold}) in one dimension. To do this in practice,
three main steps should be executed if one wishes to use methods
that partition the search region and try to obtain and to use a
Lipschitz/H\"{o}lder information   (see
\cite{pinter,dividethebest,Sergeyev_Kvasov_book} for a general
description of this kind of methods, their convergence properties,
etc.). First, in order to realize the passage from the
multi-dimensional problem to the one-dimensional one, computable
approximations to the Peano curve should be employed in the
numerical algorithms. Second, the H\"{o}lder constant $H$ from
(\ref{hold}) should be estimated. Finally, a partition strategy of
the search region should be chosen.

We have already seen above that in the Lipschitz global optimization
there exist at least four ways to obtain estimates of the Lipschitz
constant. When we move to the H\"{o}lder global optimization the
situation is different. In the literature (see
\cite{jaumard,daniya,Lera&Sergeyev2010,Lera&Sergeyev2010b,Lera&Sergeyev2013,
Sergeyev_Lera_book}), there exist methods that use   strategies
  a), b), and c) discussed above. However,
inventing for the H\"{o}lder global optimization a strategy similar
to the technique d) was an open problem since 1993 when the article
\cite{jones} showing how to use simultaneously several estimates of
$L$ has been published. In this paper, we close this gap and propose
an algorithm that uses several estimates of the H\"{o}lder constant
at each iteration employing space-filling curves, central point
based partition strategies, and H\"{o}lderian minorants.

The rest of the paper has the following structure. In Section~2,
difficulties regarding the usage of the strategy d) in the
H\"{o}lderian framework and the proposed solution are presented. In
Section~3, a new algorithm for solving the problem (\ref{p}),
(\ref{lip}) and its convergence properties are described. The new
method uses numerical approximations to Peano space-filling curves
and the scheme of representation of intervals with H\"{o}lderian
minorants from Section~2.  Section~4 presents results of numerical
experiments that compare the new method with its competitors  on 800
test functions randomly generated by the GKLS-generator from
\cite{gaviano}. Finally, Section~5 contains a brief conclusion.

\section{Lipschitz and H\"{o}lder minorants in one dimension}
\setcounter{equation}{0}

\begin{figure}[t]
\centering
\begin{minipage}{.450\linewidth}
\centerline{\resizebox{\textwidth}{!}{\includegraphics{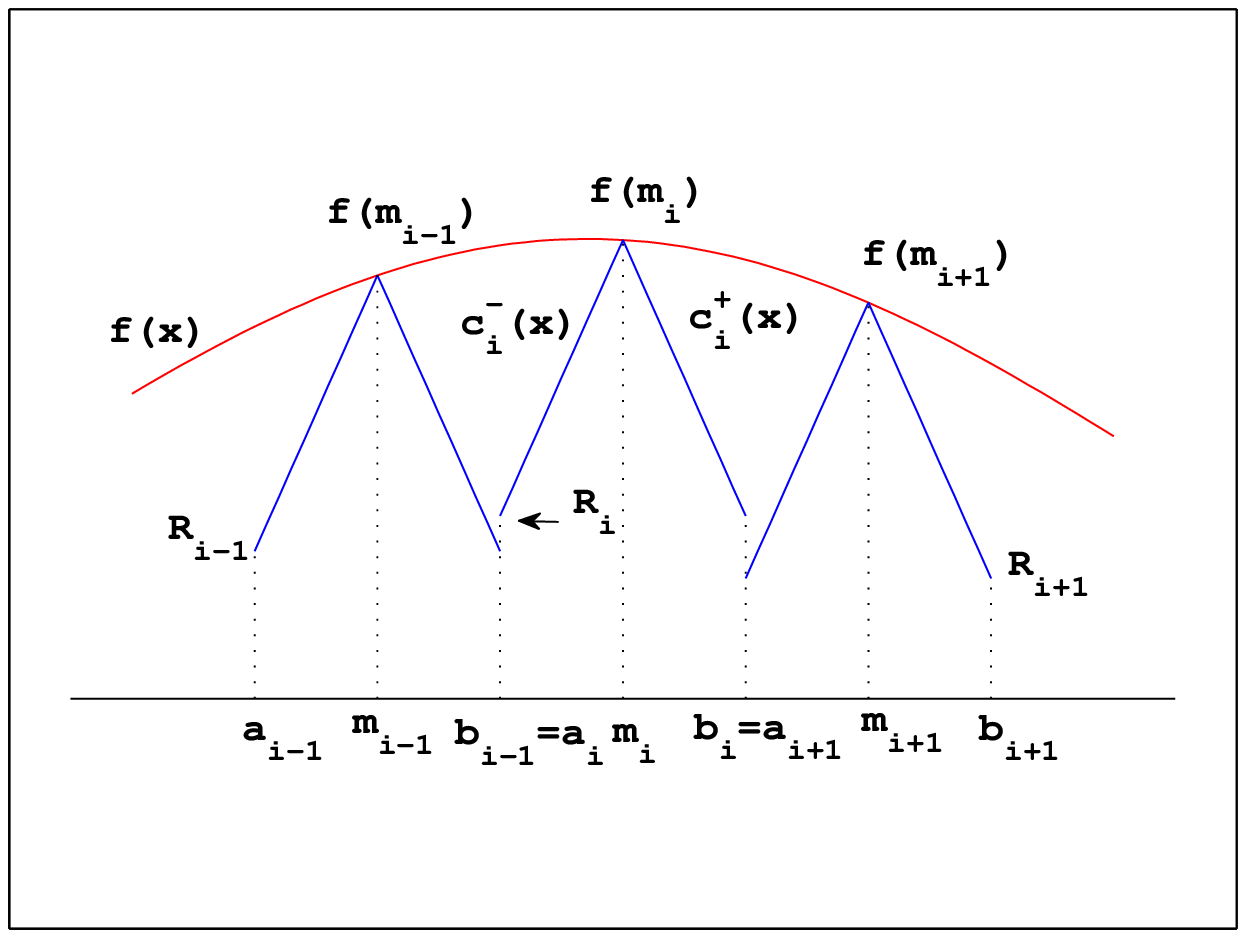}}}
\end{minipage} \hspace{0.1cm}
\begin{minipage}{.450\linewidth}
\centerline{\resizebox{\textwidth}{!}{\includegraphics{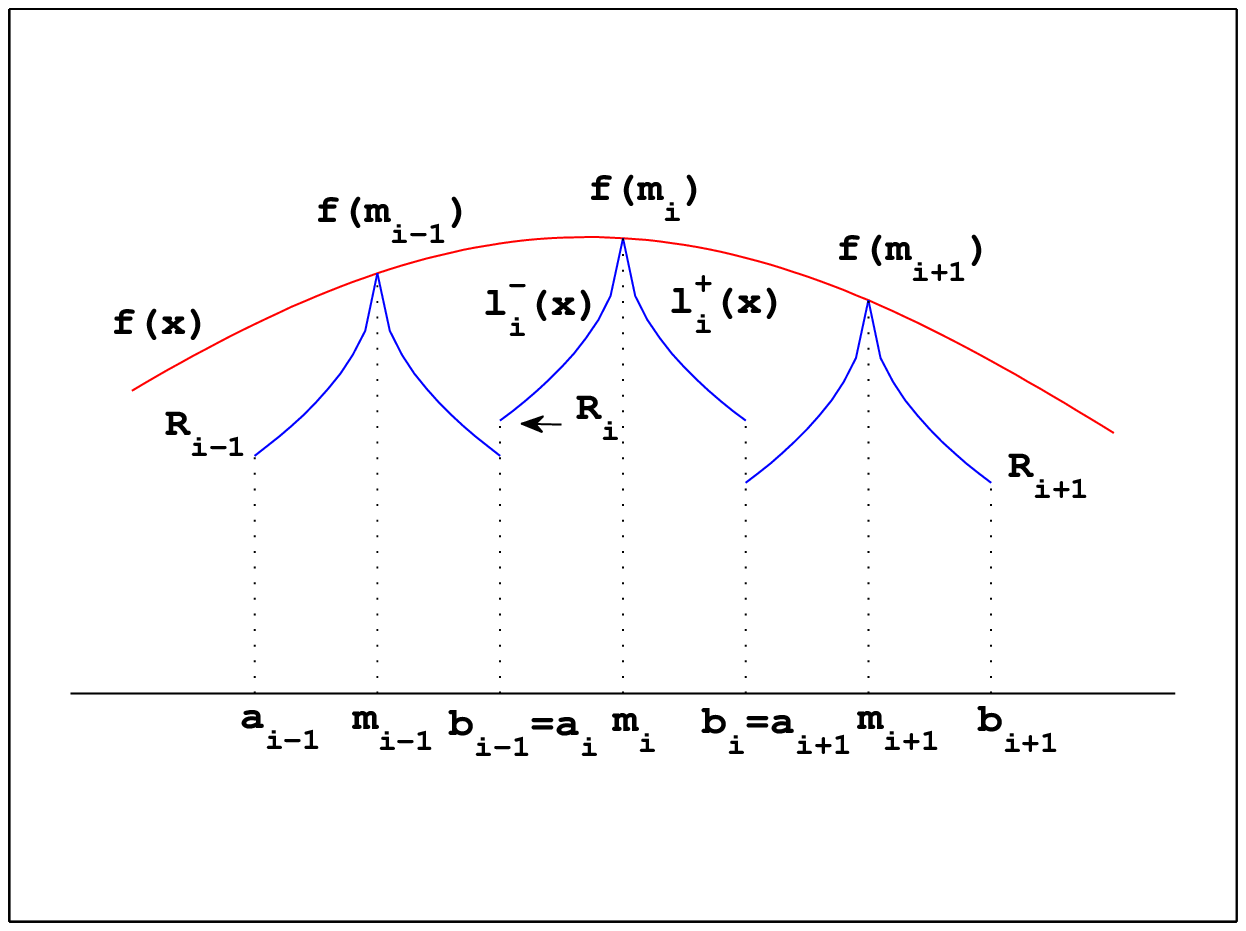}}}
\end{minipage}
\caption{\em Lipschitz (left) and H\"{o}lder (right) support
functions. } \label{fig.1}
\end{figure}

Let us consider a one-dimensional function $f(x)$ satisfying the
Lipschitz condition with a constant $L$ over the interval $[0,1]$.
Then it follows from the Lipschitz condition that
 \begin{equation}
\label{sup_Lip_1}
 f(x) \ge C^k(x) =   c_i(x),\,\,\, x\in [a_i,b_i], \,\, 1 \le i \le k,
\end{equation}
\begin{equation}
 \label{sup_Lip_2}
c_i(x) = \left\{ \begin{array}{l}
    c_i^-(x) = f(m_i) -L (m_i-x), \hspace{1.5cm}  x\in [a_i,m_i], \\
    c_i^+(x) =  f(m_i) -L (x-m_i), \hspace{1.5cm}  x\in [m_i, b_i],
    \end{array} \right.
\end{equation}
where   $C^k(x)$ is (see Fig.~\ref{fig.1}, left) a piece-wise
linear discontinuous minorant (called often also \textit{support
function}) for $f(x)$ over each subinterval $d_i=[a_i,b_i]$, $1
\le i \le k,$ and
\begin{equation}\label{m_i}
 m_i=(a_i+b_i)/2.
\end{equation}
 The values $R_i$, $1 \le i \le k,$ called
\textit{characteristics} and being lower bounds for the function
$f(x)$ over each interval $d_i$, $1 \le i \le k,$ can easily be
calculated. In fact, if we suppose that an overestimate $L_1 \ge L$
of the Lipschitz constant $L$ is given, then it follows
\begin{equation} \label{lowb_Lip}
 R_i= R_i(L_1) = \min_{x \in [a_i,b_i]} c_i(x) = f(m_i) - 0.5L_1 (b_i-a_i).
\end{equation}

However, in order to solve the multidimensional problem (\ref{p}),
(\ref{lip}) by using space-filling curves, instead of working with
Lipschitz functions, we should focus our attention on the
one-dimensional H\"{o}lderian function $f(x)$ from (\ref{fun}). It
follows from (\ref{hold}) that, for all $x, z\in [0,1]$ we have
\begin{equation} \label{low}
 f(x) \geq f(z)-H |x-z|^{1/N}.
\end{equation}
If a point $z\in [0,1]$ is fixed, then the function
$$ Z(x)=f(z)-H |x-z|^{1/N} $$
is a minorant for $f(x)$ over $[0,1]$. Then, analogously to
(\ref{sup_Lip_1})--(\ref{m_i}), we obtain that the function
\begin{equation} \label{supp}
Z^k(x) =    l_i(x),\,\,\, x\in [a_i,b_i], \,\, 1 \le i \le k,
\end{equation}
\begin{equation} \label{li}
l_i(x) = \left\{ \begin{array}{l}
    l_i^-(x) = f(m_i) -H (m_i-x)^{1/N}, \hspace{1.5cm}  x\in [a_i,m_i], \\
    l_i^+(x) =  f(m_i) -H (x-m_i)^{1/N}, \hspace{1.5cm}  x\in [m_i, b_i],
    \end{array} \right.
\end{equation}
is a discontinuous nonlinear minorant for $f(x)$ (see
Fig.~\ref{fig.1}, right). The values $R_i$, $1 \le i \le k,$ are
lower bounds for the function $f(x)$ over each interval $d_i$, $1
\le i \le k,$ and can be  calculated as follows if an overestimate
$H_1 \ge H$ of the H\"{o}lder constant $H$ is given
\begin{equation} \label{lowb}
 R_i= R_i(H_1) = \min_{x \in [a_i,b_i]} l_i(x) = f(m_i) -H_1 |(b_i-a_i)/2|^{1/N}.
\end{equation}

\begin{figure}[t]
\centerline{\epsfxsize=14.5cm\mbox{\epsfbox{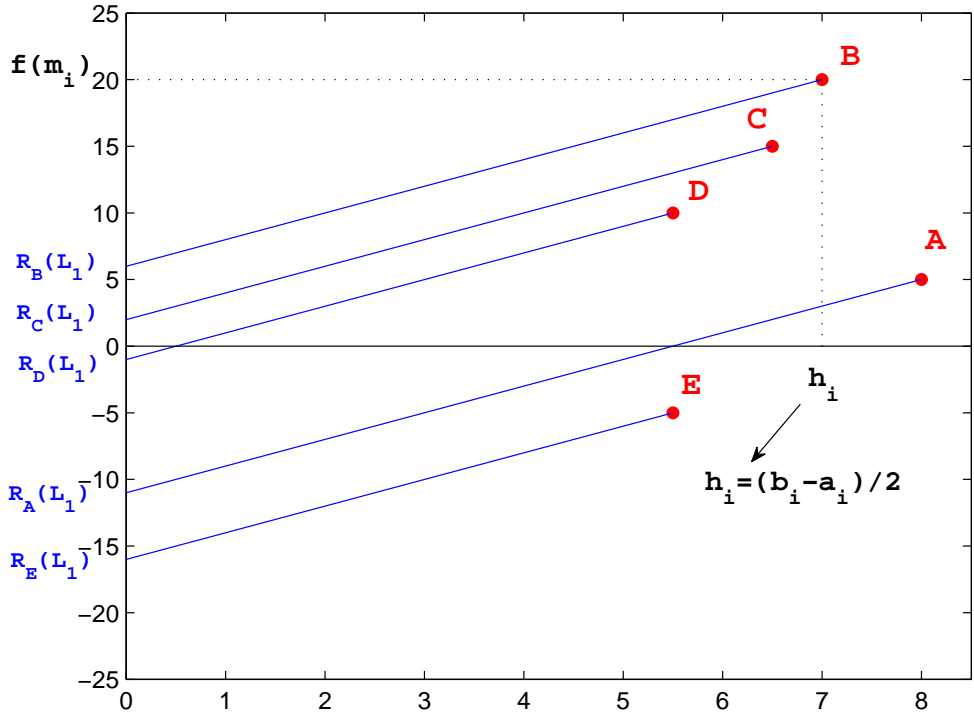}}}
\caption{\em Representation of intervals in the Euclidean metric in
the framework of the DIRECT algorithm.} \label{direct}
\end{figure}

As it was discussed above, the DIRECT algorithm works simultaneously
with several estimates of the Lipschitz constant at each iteration.
One of the key features that allow it to do this is a smart
representation of intervals $[a_i,b_i]$, $1 \le i \le k,$ in the
two-dimensional diagram shown in Fig.~\ref{direct}. In this Figure,
we have   intervals $d_A$, $d_B$, $d_C$, $d_D$, and $d_E$ that are
represented by the dots $A$, $B$, $C$, $D$, and $E$, respectively.
The horizontal coordinate of each dot is equal to $0.5(b_i-a_i)$,
i.e., half of the length of the respective interval $[a_i,b_i]$, and
the vertical coordinate is equal to $f(m_i)$ where $m_i$ is from
(\ref{m_i}) (see, e.g., the dot $B$ and its coordinates). Let us
consider now the intersection of the vertical coordinate axis with
the line having the slope $L_1$ and passing through each dot
representing subintervals in the diagram shown in Fig.~\ref{direct}.
It is possible to see that this intersection gives us exactly the
characteristic $R_i=R_i(L_1)$ from~(\ref{lowb_Lip}), i.e., the lower
bound for $f(x)$ over the corresponding subinterval $[a_i,b_i]$ if
$L_1 \ge L$ (note that the points on the vertical axis ($d_i=0$) do
not represent any subinterval).

It can immediately be  seen that the diagram allows one to determine
very easily an interval with the minimal characteristic (in
Fig.~\ref{direct} this interval is represented by the dot $E$, its
characteristic is $R_E(L_1)$). In the Lipschitz global optimization
(see, e.g.,
\cite{pinter,dividethebest,Sergeyev_Kvasov_book,Sergeyev_Lera_book,yaro}),
the operation of finding an interval with the minimal characteristic
(that can be calculated in different ways in different algorithms)
is an important one. It is executed at each iteration, the found
interval is  then subdivided and $f(x)$ is evaluated at new points
belonging to this interval. Moreover, since we do not know the exact
value of the real Lipschitz constant $L$, the scheme presented in
Fig.~\ref{direct} allows one (see \cite{jones}) to take into
consideration all possible values of $L$ from zero to
infinity\footnote{In fact, it is easy to see that if we connect the
points $A$ and $E$ by a line and indicate its slope by $L^{max}$
then for all constants $L \ge L^{max}$ the interval  $d_A$ will have
the minimal characteristic and, therefore it should be subdivided.
Analogously, for all constants $L \le L^{max}$ the interval  $d_E$
will have the minimal characteristic. Then, if we subdivide both
intervals, $d_A$ and $d_E$, then the interval corresponding to the
real constant $L$ will be subdivided even though we are not know
this value $L$. Thus, the diagram really allows one to consider all
possible values of the Lipschitz constant from zero to infinity and
to choose a small group of intervals for the subsequent subdivision
ensuring that the interval corresponding to the correct Lipschitz
constant belongs to this group. } and to choose a set of promising
intervals for partitioning (this issue will be discussed in detail
later).

\begin{figure}[t]
\centerline{\epsfxsize=14.5cm\mbox{\epsfbox{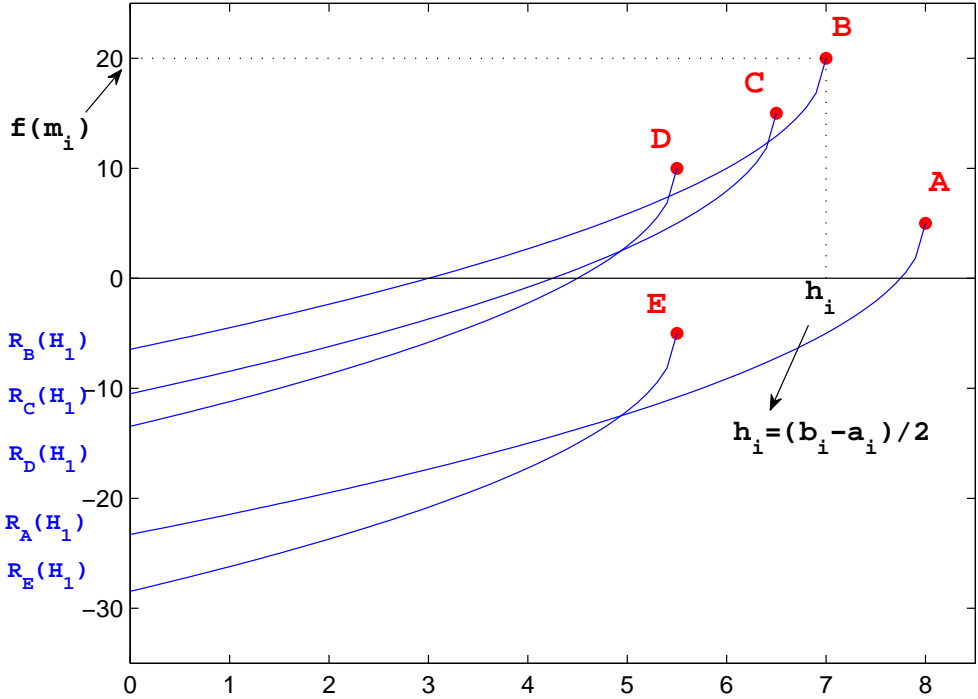}}}
\caption{\em Due to numerous possible intersections of nonlinear
H\"{o}lderian minorants, the attempt to use the representation of
intervals analogous to that  shown in Fig.~\ref{direct} does not
give the desired effect.} \label{fig.2}
\end{figure}

Let us try  now to construct a similar diagram in the framework of
the H\"{o}lder optimization where the nonlinear support functions
$l_i(x)$ from (\ref{li}) shown in Fig.~\ref{fig.1}, right, are built
over each subinterval.  In Fig.~\ref{fig.2}, intervals $d_A$, $d_B$,
$d_C$, $d_D$, and $d_E$ are again represented by dots $A$, $B$, $C$,
$D$, and $E$, respectively. If we take an estimate $H_1$ of the
H\"{o}lder constant  $H$, then characteristic $R_B(H_1)$ of the
interval $d_B$ represented by the dot $B$ is given by the vertical
coordinate of the intersection point of the auxiliary curve
(\ref{li}) passed through the point $B$ and the vertical coordinate
axis.

We can see in Fig.~\ref{fig.2} that the curves constructed using the
estimate  $H_1$ and representing the nonlinear  support function
(\ref{li}) can intersect each other in different ways. The procedure
of the selection of subintervals for producing new trials
(\textit{trial} is an evaluation of $f(x)$ at a point $x$ that is
called \textit{trial point}) is based on estimates of the lower
bounds of the objective function over current subboxes. Subintervals
for the further subdivision are chosen taking into account
\textit{all} possible values of the H\"{o}lder constant $H$. Due to
numerous possible intersections of the curves at the representation
shown in Fig.~\ref{fig.2}, it becomes unclear how to select, for all
possible values of $H$, the set of promising intervals that must be
partitioned at each iteration~$k$.

This difficulty that seemed unsolvable since 1993 did not allow
people to construct global optimization algorithms working in the
framework of the H\"{o}lder global optimization with all possible
H\"{o}lder constants. This problem is solved in the next section
that proposes an algorithm that uses, on the one hand, Peano
space-filling curves to reduce the multi-dimensional problem
(\ref{p}), (\ref{lip}) to one dimension  and, on the other hand,
multiple estimates of the  H\"{o}lder (Lipschitz) constant. The
algorithm for solving one-dimensional H\"{o}lder global optimization
problems in the spirit of the DIRECT method is a part of it.

\section{The new algorithm}
\setcounter{equation}{0}

 In order to construct a procedure allowing one to select,
at each iteration $k$ of the algorithm, a set of intervals to be
partitioned with respect to all possible values of the constant
$H$  we proceed as follows. First, we introduce a new graphical
representation of  subintervals $d_i$ by using instead of
Euclidean metric  the H\"{o}lderian one. Namely, we propose to
represent each interval $d_i=[a_i,b_i] \in\{ D^k\}$  by a dot
$P_i$ with coordinates $(h_i,F_i)$, where $D^k$ is the current
partition of the one-dimensional search interval during the $k$th
iteration and coordinates of the point $P_i$  are calculated as
follows
\begin{equation}\label{xdot}
 h_i=|(b_i-a_i)/2|^{1/N},
 \end{equation}
\begin{equation}\label{dot}
 F_i=f(m_i),
\end{equation}
where $m_i$ is from (\ref{m_i}), i.e., $m_i$ is the central point
of the interval $d_i=[a_i,b_i]$.

\begin{figure}[t]
\centerline{\epsfxsize=14.5cm\mbox{\epsfbox{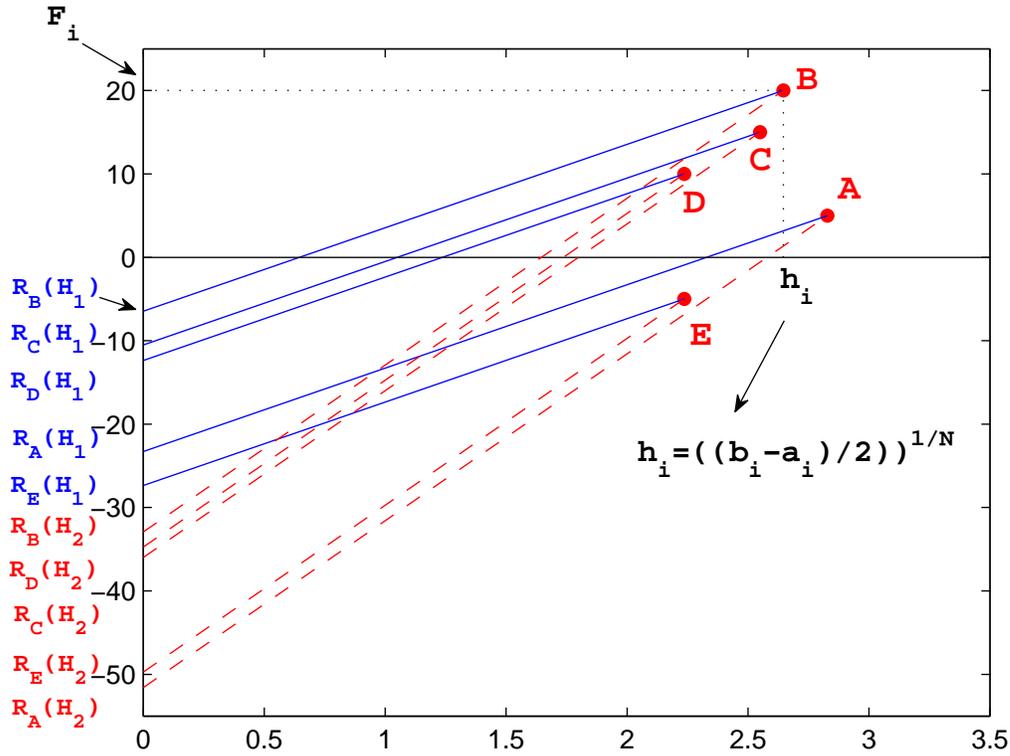}}}
\caption{\em Representation of intervals with the H\"{o}lderian
metric.} \label{fig.3}
\end{figure}

The introduction of the H\"{o}lderian metric allows us (see
Fig.~\ref{fig.3}) to avoid the non-linearity and the intersection of
minorants giving as a result a diagram similar to that   from
Fig.~\ref{direct}. In Fig.~\ref{fig.3} we can see the representation
of the same intervals considered in Fig.~\ref{fig.2}. Notice that in
Fig.~\ref{fig.3} the values in the horizontal axis are calculated in
the H\"{o}lderian metric, while the vertical axis values coincide
with those of Fig.~\ref{fig.2}. In this new representation, the
intersection of the line with the slope $H_1$ passing through any
dot representing an interval in the diagram of Fig.~\ref{fig.3} and
the vertical coordinate axis gives us exactly the characteristic
(\ref{lowb}) of the corresponding interval.

We can proceed now with the development of the new one-dimensional
global optimization method following the spirit of the DIRECT
algorithm and keeping in mind that we deal with the H\"{o}lderian
metric.  Each subinterval $d_i$ of a current partition $D^k$ is
characterized by a lower bound of the objective function over $d_i$.
An interval $d_i$ is selected for a further partitioning  if for
some value $\tilde{H}>0$ of the possible H\"{o}lder constant it has
the smallest lower bound for $f(x)$ with respect to the other
intervals present at this iteration. By changing the value of
$\tilde{H}$ from zero to infinity, at each iteration $k$, we select
a set of intervals that will be partitioned.

Let us consider how this set of intervals is chosen  during each
iteration $k$. The following definitions state a \textit{relation
of domination} between every two subintervals of the current
partition of~$D$.
\begin{defin}\label{def1}
 Given an estimate $\tilde{H}>0$ of the H\"{o}lder constant $H$ from (\ref{hold}), an interval
$d_i\in \{ D^k\}$ dominates the interval $d_j\in \{ D^k\}$ with
respect to $\tilde{H}$ if
$$ R_i(\tilde{H}) < R_j(\tilde{H}).  $$
\end{defin}
\begin{defin}\label{def2}
An interval $d_t\in \{ D^k\}$ is said to be nondominated with
respect to $\tilde{H}>0$ if for the chosen value $\tilde{H}$ there
is no other interval in $\{ D^k\}$ which dominates~$d_t$.
\end{defin}

In order to be sure to subdivide the nondominated interval
corresponding to the real constant $H$, we can select   the set of
nondominated intervals with respect to all possible estimates
$0<\tilde{H}<\infty$. By using the new graphical representation
shown in Fig.~\ref{fig.3} it is easy to determine whether an
interval dominates, with respect to an estimate of the H\"{o}lder
constant $H$, some other interval from the partition $\{ D^k\}$.
For example, in Fig.~\ref{fig.3} we can see that for the estimate
$H_1$ we have
$$ R_D(H_1)<R_C(H_1)<R_B(H_1), $$
so, the interval $d_D$ dominates both intervals $d_C$ and $d_B$ with
respect to $H_1$. Furthermore, $ R_A(H_1)<R_D(H_1)$, i.e., the
interval $d_A$ dominates the interval $d_D$. Finally, the
characteristic $R_E(H_1)$ is the smallest one, and  the interval
$d_E$ dominates all others intervals with respect to $H_1$, i.e., it
is nondominated  with respect to $H_1$ (see Fig.~\ref{fig.3}).

 If a
different value $H_2>H_1$ of the H\"{o}lder constant is considered
(see Fig.~\ref{fig.3}), the interval $d_D$ still dominates the
interval $d_B$ with respect to  $H_2$, because $R_D(H_2)<R_B(H_2)$,
but $d_D$ is dominated by the interval $d_C$, since
$R_D(H_2)>R_C(H_2)$. Moreover, we have that
$R_A(H_2)<R_E(H_2)<R_C(H_2)$, thus, for the chosen estimate $H_2$
the unique nondominated interval is $d_A$. As we can see from this
example, some intervals ($d_B$, $d_C$, and $d_D$ in
Fig.~\ref{fig.3}) are always dominated by other intervals,
independently on the chosen estimate of the H\"{o}lder constant.
Other intervals ($d_E$ and $d_A$ in Fig.~\ref{fig.3}) can be
nondominated for one value and dominated for another one.  The
following definition will be very useful hereinafter.

\begin{defin}\label{def3}
 An interval $d_t\in \{ D^k\}$ is called nondominated if there exists an estimate $0<\tilde{H}<\infty$ of
the H\"{o}lder constant $H$ such that $d_t$ is nondominated with respect to $ \tilde{H}$.
\end{defin}

\begin{figure}[t]
\centerline{\epsfxsize=14.5cm\mbox{\epsfbox{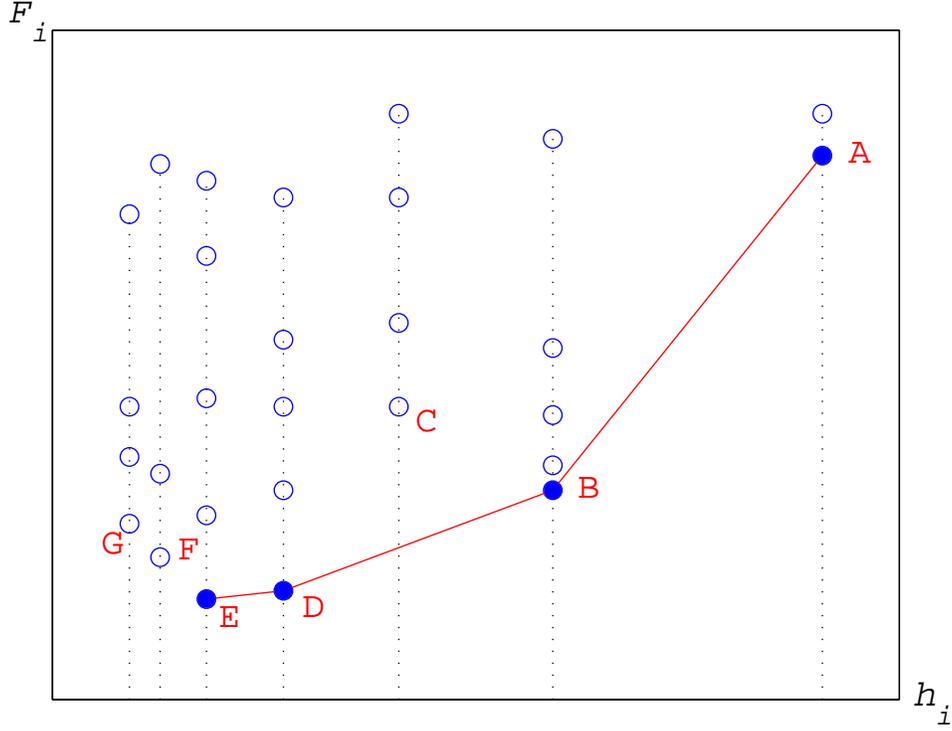}}}
\caption{\em The nondominated intervals $d_A$, $d_B$, $d_D$, and
$d_E$  are represented by   dots $A$, $B$,
 $D$, and  $E$.} \label{fig.4}
\end{figure}

In other words, nondominated intervals are intervals over which the
objective function $f(x)$ has the smallest lower bound from
(\ref{lowb}) for some particular estimate of the H\"{o}lder constant
$H$. Note that in the two-dimensional diagram $(h_i,F_i)$, where
$h_i$ and $F_i$ are from (\ref{xdot}), (\ref{dot}), the nondominated
intervals are located at the bottom of each group of points with the
same horizontal coordinate. For example in Fig.~\ref{fig.4} these
points are designated as $A$, $B$, $C$, $D$,  $E$, $F$,  and $G$.
Not all of these intervals are nondominated: in fact, in
Figure~\ref{fig.4} the interval $d_C$ is dominated either by the
interval $d_B$ (for example, with respect to $H_1 \geq H_{BD}$,
where $H_{BD}$ corresponds to the slope of the line passed through
the points $B$ and $D$), or by the interval $d_D$ (with respect to
$H_2 < H_{BD}$). The interval $d_F$ is dominated by $d_D$ and $d_E$,
with respect to any positive estimate of the constant $H$. The
interval $d_G$ is dominated by $d_F$.  In Figure \ref{fig.4}, dots
$A$, $B$,
 $D$, and  $E$  represent nondominated
intervals. The following theorem allow us to easily identify the
set of nondominated intervals.
\begin{theorem}\label{th0}
Let each interval $d_i=[a_i,b_i] \in \{ D^k\}$ be represented by a dot with horizontal
coordinate $h_i$ and vertical coordinate $F_i$ defined in (\ref{xdot}), (\ref{dot}).
Then, intervals that are nondominated in the sense of Definition \ref{def3} are located
on the lower-right convex hull of the set of dots representing the intervals.
\end{theorem}

\textbf{Proof.} The proof of the Theorem \ref{th0} is analogous to
the proof of Theorem 2.2 from \cite{Sergeyev_Kvasov(2006)}. \hfill
$\Box$

In practice,  nondominated intervals can be found by applying
algorithms for identifying the convex hull of the dots (see, e.g.,
the algorithm called Jarvis march, or gift wrapping, see
\cite{preparata}).

We describe now the partition strategy adopted by the new
algorithm for dividing subintervals in order to produce new trial
points. When, at the generic iteration $k$, we identify the set of
nondominated intervals, we proceed with the subdivision of each of
these intervals only if a significant improvement on the function
values with respect to the current minimal value $f_{min}(k)$ is
expected. Once an interval $d_t\in \{ D^k\} $ becomes
nondominated, it can be subdivided only if the following condition
is satisfied:
\begin{equation} \label{cond}
 R_t(\tilde{H}) \leq f_{min}(k) -\xi,
\end{equation}
where the lower bound $R_t=R_t(\tilde{H})$ is from (\ref{lowb}).
This condition prevents the algorithm from subdividing already
well-explored small subintervals.

Let us suppose now that at the current iteration $k$ of the new
algorithm a subinterval $d_t=[a_t,b_t]$, represented in the
two-dimensional diagram of Fig.~\ref{fig.4} by the dot $(h_t,F_t)$
from (\ref{xdot}), (\ref{dot}), has been chosen for partitioning.
The subdivision of this interval is performed in such a way that
three new equal subintervals of the length $(b_t-a_t)/3$ are
created, i.e.,
\begin{equation} \label{partit1}
 [a_t,b_t]=[a_t,p_t]\cup [p_t, q_t] \cup [q_t, b_t],
\end{equation}
\begin{equation} \label{partit2}
 p_t=a_t+(b_t-a_t)/3, \ \ q_t=b_t-(b_t-a_t)/3.
\end{equation}
The interval $d_t$ is removed from the two-dimensional diagram
representing the current partition $\{ D^k\}$ of the search
interval, and the three newly generated subintervals are introduced
both into $\{ D^k\} $ and the diagram. Finally, two new trials
$f(c_1)$ and $f(c_2)$ are performed at the   points $c_1$ and $c_2$
of the intervals $ [a_t,p_t]$ and $[q_t,b_t]$, respectively, i.e.,
\begin{equation} \label{trials}
 c_1=\frac{a_t+p_t}{2}, \hspace{2cm} c_2=\frac{q_t+b_t}{2}.
\end{equation}
The central interval $[p_t,q_t]$ inherits the point
$\frac{a_t+b_t}{2}=\frac{p_t+q_t}{2}$ at which the objective
function has been evaluated when the original interval
$d_t=[a_t,b_t]$ has been created.

Until now we have described the strategies assuming to work with a
function in one dimension. As was stated above, Strongin has shown
that multidimensional optimization problems  can be solved by using
modified  algorithms proposed for minimizing functions in one
dimension, and therefore in order to solve the global optimization
problem in $N$ dimensions (\ref{p}), (\ref{lip}) we can use the
above developed one-dimensional global optimization method together
with the space-filling curves. For an effective use of the Peano
curve in our algorithm we need computable approximations of the
curve (see \cite{Sergeyev_Lera_book,yaro} for a detailed discussion
and a code allowing one to implement such approximations).
Hereinafter we denote by $p_M(x)$ the approximation of level $M$ of
the Peano curve. In Fig.~\ref{fig.1a} we can see examples of Peano
curve approximations of the levels $M=2,3,4,5$ in dimension $N=2$.

Suppose now that a global optimization method uses an approximation
$p_M(x)$ of the Peano curve  to solve the multidimensional problem
and provides a lower bound $U_M^*$ for the corresponding
one-dimensional function $f(x)$. Then the value $U_M^*$ will be a
lower bound for the function $F(y)$ in dimension $N$ only along the
curve $p_M(x)$. The following theorem establishes a lower bound for
the function $F(y)$ over the entire multidimensional search region
$[a,b]$ given the value $U_M^*$ .
\begin{theorem} \label{th1}
 Let $U^*_M$ be a lower bound along the space-filling curve $p_M(x)$ for a
 multidimensional function $F(y)$, $y\in [a,b]\subset R^N$, satisfying Lipschitz
 condition with constant $L$, i.e.,
 $$ U^*_M \leq F(p_M(x)), \hspace{1cm} x\in[0,1].  $$
 Then the value
 $$ U=U^*_M - 2^{-(M+1)}L\sqrt N   $$
 is a lower bound for $F(y)$ over the entire region $[a,b]$.
\end{theorem}

\textbf{Proof.} See \cite{Lera&Sergeyev2010b} or the recent
monograph \cite{Sergeyev_Lera_book}. \hfill $\Box$

By using the space-filling curves we are able to work with a
one-dimensional function in the interval $[0,1] \subset R^1$. The
level $M$ of the approximation of the Peano curve  $p_M(x)$, is
crucial for a good performance of the method. If $M$ is too small,
the domain in $N$ dimensions may not be well ``filled'' and we risk
losing the global solution. When $M$ grows, the reduced function in
one dimension becomes more and more oscillating  and the number of
local minima increases when  $N$ increases (see
\cite{Lera&Sergeyev2010b} for a detailed discussion on this topic).
Then, due to the facts that we are in $[0,1]$ and that we use the
metric of H\"{o}lder, it happens that the width of the nondominated
interval $d_t\in \{ D^k\}$ to be partitioned at a generic iteration
$k$  can become very small. When the dimension $N$ increases, the
width of the subintervals can reach the computer precision. In order
to avoid this situation another condition in addition to
(\ref{cond}) is required. Namely, when an interval $d_t=[a_t,b_t]\in
\{ D^k\} $ becomes nondominated, it can be subdivided only if the
following condition is satisfied:
\begin{equation} \label{cond1}
 b_t-a_t > \eta,
\end{equation}
where $\eta$ is a parameter of the method.

Now, let us present  the new algorithm  called $MGAS$
($M$ultidimensional $G$lobal optimization $A$lgorithm working with a
$S$et of estimates of the H\"{o}lder constant).

\vskip 24pt
{\bf  Algorithm MGAS}
\begin{description}
\item[Step 0.] (Initialization). Set the current iteration number $k:=1$.

Split the initial interval $D=[0,1]$ in three equal parts and set
$x^1=1/6$, $x^2=1/2$, $x^3=5/6$ and compute the values of the
function $z^j=f(x^j)=F(p_M(x^j))$, $j=1,2,3$, where $p_M(x)$ is the
$M$-approximation of the Peano curve.

Set the current partition of  the search interval $D^1=\{ [0,1/3],
[1/3,2/3], [2/3,1] \}$.

Set the current number of intervals $J=3$ and the current number of
trials $T=3$.

Set $f_{min}(1)=\min \{z^1, z^2, z^3 \}$,  and $x_{min}(1)=argmin
\{f(x^j) : j=1,2,3 \}$.

After executing $k$ iterations, the iteration $k+1$ consists of
the following steps.

\item[Step 1.] (Nondominated intervals) Identify both the set $S^k$, $S^k\subset D^k$, of nondominated intervals,
according to Definition \ref{def3}, that satisfy conditions
(\ref{cond}) and (\ref{cond1}), and the corresponding set $I^k$ of
their indices. $D^k$ denotes the partition of the search interval
$D=[0,1]$ at iteration $k$.
\item[Step 2.] (Subdivision of nondominated intervals) Set $D^{k+1}=D^k$, and perform the following Steps 2.1-2.3:
\begin{flushleft}
\item[{\bf Step 2.1}] (Interval selection). Select a new interval $d_t=[a_t,b_t]$ from $S^k$  such that
$$ t=arg \max_{j\in I^{k}}\{ b_j-a_j \}. $$
 \item[{\bf Step 2.2}] (Subdivision and sampling). Subdivide interval $d_t$ in three  new equal subintervals,
named $d_{t1}$, $d_{t2}$, $d_{t3}$, of the length $(b_t-a_t)/3$,
following (\ref{partit1}), (\ref{partit2}),
and produce two new trial points accordingly to (\ref{trials}). \\
Eliminate the interval $d_t$ from $D^{k+1}$, i.e., set $D^{k+1}=D^{k+1}\setminus \{ d_t\} $, and update
$D^{k+1}$ with the insertion of the  three new intervals, i.e.,
$$D^{k+1}= D^{k+1}\cup \{ d_{t1}\} \cup \{ d_{t2}\} \cup \{ d_{t3}\}. $$
Increase both the current number of intervals $J=J+2$, and
 the current number of trials $T=T+2$. \\
Update the current record $f_{min}$ and the current record point $x_{min}$ if
necessary.
\item[{\bf Step 2.3}] (Next interval). Eliminate the interval $d_t$ from $S^k$, i.e., set $S^k=S^k\setminus \{ d_t\}$
and $I^k=I^k\setminus \{ t\}$.

If $S^k \neq \emptyset$, then go to Step 2.1. Otherwise go to Step 3.
\end{flushleft}
\item[Step 3.] (End of the current iteration). Increase the iteration counter $k=k+1$. Go to Step 1 and start
the next iteration.
\end{description}

Different stopping criteria can be used in  the algorithm introduced
above. One of them will be introduced in  the next section
presenting numerical experiments.

 We proceed now to the study of convergence properties of the algorithm. Theorem
 \ref{th1} linking the multidimensional global optimization problem (\ref{p}), (\ref{lip})
 to the one-dimensional problem (\ref{fun}), (\ref{hold}) allows us to concentrate our
attention on the convergence properties of the one-dimensional
method.
 We shall study properties of an infinite sequence $\{x^{j(k)}\}$ of trial points
 generated by the algorithm MGAS when we suppose that the number of iteration
 $k$ goes to infinity (i.e., in this case the algorithm does not stop).
 The following theorem establishes the so-called \textit{everywhere dense} convergence of the
 method, i.e., convergence of the infinite sequence of trial points to any point
 of the one-dimensional search domain.
\begin{theorem} \label{th2}
If $\eta=0$ in (\ref{cond1}), then for any point $x\in D=[0,1]$ and
any $\delta >0$ there exist an iteration number
 $k(\delta)\geq 1$ and a point $x'\in \{x^{j(k)}\}$, $k>k(\delta)$, such that
 $ |x-x'|< \delta $.
\end{theorem}
 {\em Proof}.
 The interval partition scheme (\ref{partit1}), (\ref{partit2}) used for each subdivision
 of  intervals produces three new subintervals of the length equal to a third of the
 length of the subdivided interval. Since $\eta=0$, to prove the Theorem  it is sufficient
 to prove that for a fixed value $\delta >0$, after a finite number of iterations $k(\delta)$, the largest subinterval
 of the current partition $\{D^{k(\delta)}\}$ of the domain $D$ will have the length
 smaller than~$\delta$. In such a case, in the $\delta$-neighborhood of any point of $D$
 there will exist at least one trial point generated by the
 algorithm. $\hfill \Box$

To see this, let us  fix an iteration number $k'$ and  consider
the group of the largest intervals of the partition $\{D^{k'}\}$
having the horizontal coordinate $h_{max}$ (in the diagram of
Fig.~\ref{fig.4} this group consists of two points: the dot $A$
and the dot above it). As can be seen from the scheme of the
algorithm MGAS, for any $k'\geq 1$ this group is always taken into
account when nondominated intervals are looked for. In particular,
an interval $d_t \in \{D^{k'}\}$ from this group, having the
smallest value $F_t$, must be partitioned at each iteration of the
algorithm. This happens because there always exists a sufficiently
large estimate $H_\infty$ of the H\"{o}lder constant $H$ for the
function $f(x)$ such that the interval $d_t$ is the nondominated
interval with respect to $H_\infty$ and condition (\ref{cond}) is
satisfied for the lower bound $R_t(H_\infty)$.

Three new subintervals having the length equal to a third of the
length of $d_t$ are then inserted into the group with a horizontal
coordinate $h_j < h_{max}$. Since each group contains a finite
number of intervals, after a sufficiently large number of iterations
all the intervals of the group $h_{max}$ will be divided and the
group will become empty. As a consequence, the group of the largest
intervals will now be identified by $h_j$, where the difference
$h_{max}-h_j>0$ is finite. The same procedure will be repeated with
this new group of the largest intervals, and the next new group,
etc.

 This means that
there exists a finite iteration number $k(\delta)$ such that after
performing $k(\delta)$ iterations of the algorithm MGAS, the
length of the largest interval of the current partition
$\{D^{k(\delta)}\}$  is smaller than $\delta$ and, therefore,  in
the $\delta$-neighborhood of any point of the search region
 there will exist at least one trial point generated by the algorithm.

\begin{figure}[!tbp]
\centerline{\epsfxsize=10cm\mbox{\epsfbox{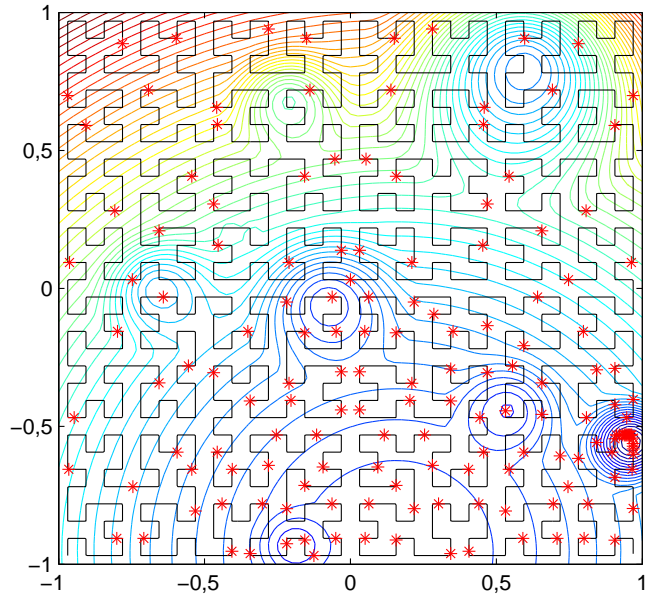}}}
\caption{\em  Trial points produced by the MGAS and Peano curve
approximation of level 5 while optimizing Function no. 6, class 1,
generated by the GKLS-generator; Table \ref{table1} describes this
and other classes of test functions used in the numerical
experiments.} \label{fig.5}
\end{figure}

In Fig.~\ref{fig.5},  an example of convergence of the sequence of
trial points generated by the algorithm MGAS in dimension $N=2$
using the approximation of the level $M=5$ to the Peano curve is
given. The zone with the high density of the trial points
corresponds to the global minimizer.

\begin{figure}[t]
\centering
\begin{minipage}{.600\linewidth}
\centerline{\resizebox{\textwidth}{!}{\includegraphics{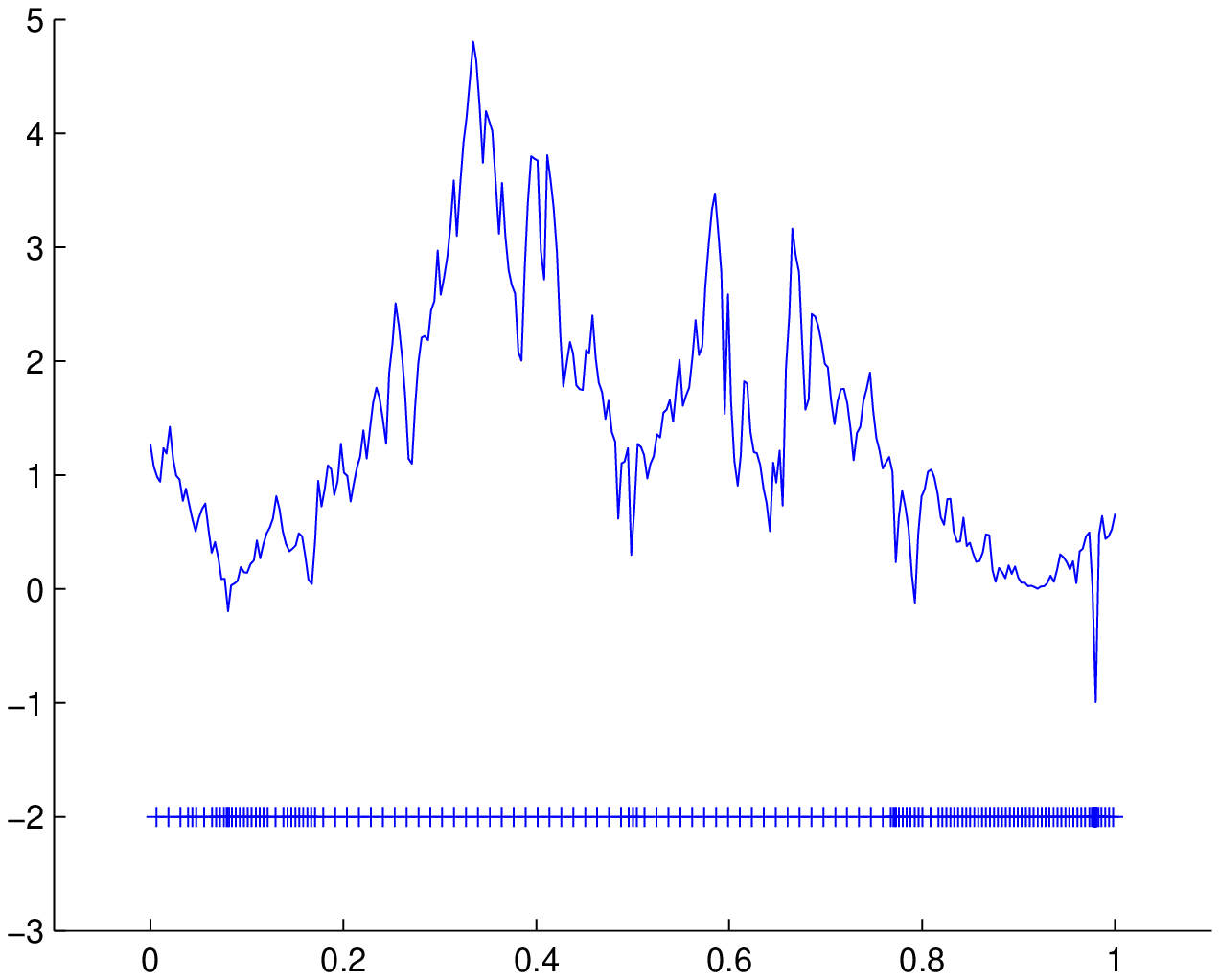}}}
\end{minipage} \hspace{0.1cm}
\begin{minipage}{.600\linewidth}
\centerline{\resizebox{\textwidth}{!}{\includegraphics{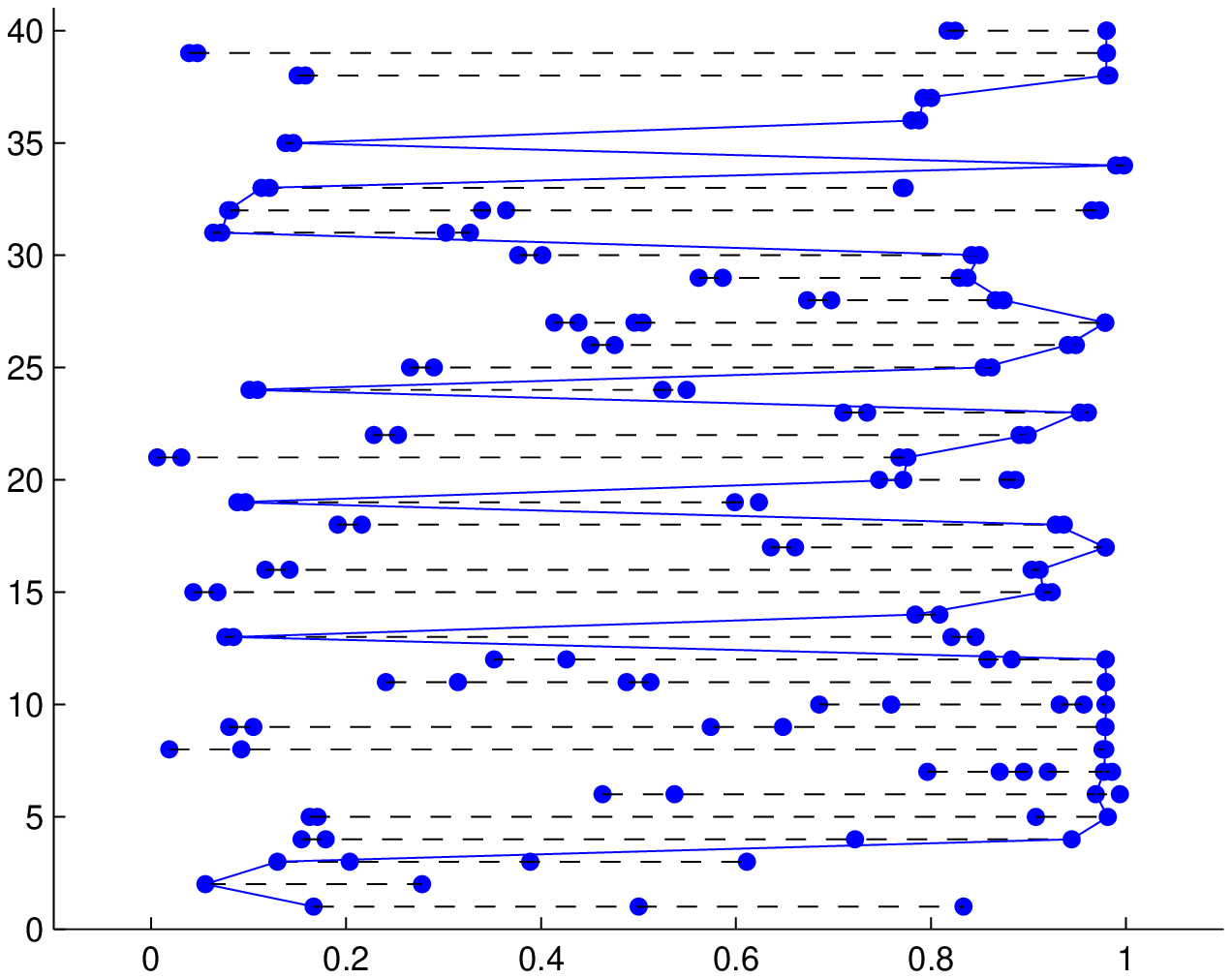}}}
\end{minipage}
\caption{\em  The one-dimensional function corresponding to the
curve shown in Fig.~\ref{fig.5} and the  respective trial points
produced by the MGAS at the interval $[0,1]$; The lower part of the
Figure shows the dynamics   of 40 iterations executed by the MGAS.}
  \label{fig.51}
\end{figure}

 Fig.~\ref{fig.51} shows how
this problem was solved in the one-dimensional space.  In the upper
part of Fig.~\ref{fig.51}, the one-dimensional function
corresponding to the curve shown in Fig.~\ref{fig.5} and the
respective trial points produced by the MGAS at the interval $[0,1]$
are presented.   The lower part of the Figure shows the dynamics
(from bottom to top) of 40 iterations executed by the MGAS. It can
be seen that each iteration contains more than one trial. The
piece-wise line connects points with the best function value
obtained during that iteration.

In order to conclude this section it should be noticed that Theorem
\ref{th2} establishes convergence conditions of infinite sequences
of trial points generated by the algorithm MGAS to any point of the
domain $[0,1]$ and therefore to the global minimum points $x^*$ of
the one-dimensional function $f(x)$. The Peano curves used for
reduction of dimensionality establish a correspondence between
subintervals of the curve and the $N$-dimensional subcubes of the
domain $[a,b] \subset R^N$. Every point on the curve approximates an
$\varepsilon$-neighborhood in $[a,b]$, i.e., the points in the
$N$-dimensional domain may be approximated differently by the points
on the curve in dependence on the mutual disposition between the
curve and the point in $[a,b]$ to be approximated (see
\cite{Sergeyev_Lera_book,yaro}). Here by approximation of a point
$y\in[a,b]$ we mean the set of points (called \textit{images}) on
the curve minimizing the Euclidean distance from $y$. It was shown
in \cite{Sergeyev_Lera_book,yaro} that the number of the images
ranges between $1$ and $2^N$. These images can be located on the
curve very far from each other despite their proximity in the
$N$-dimensional space. Thus, by using the space-filling Peano curve
$p_M(x)$, the global minimizer $y^*$ in the $N$-dimensional space
can have up to $2^N$ images on the curve, i.e., it is approximated
by $n$, $1\leq n\leq 2^N $, points $y^{*i}$ such that
$$ y^{*i} = p_M(x^{*i}), \hspace{1cm} \|y^{*i}-y^* \|\leq \varepsilon,
 \hspace{1cm} 1\leq i \leq n, $$
where $\varepsilon >0$ is defined by the space-filling curve.
Obviously, in the limiting case, when $M\rightarrow \infty$ and the
iteration number $k\rightarrow \infty$, all global minimizers will
be found. But in practice we work with a finite $M<\infty$ and
$k<\infty$, i.e., with a finite trial sequence, then to obtain an
$\varepsilon$-approximation $y^{*i}$ of the solution $y^*$ it is
sufficient to find \textit{only one} of the images $x^{*i}$ on the
curve. This effect may result in  a serious acceleration of the
search (see \cite{yaro} for a detailed discussion).

\section{Numerical experiments}
\setcounter{equation}{0}

In this section, we present results of numerical experiments
performed to compare the new algorithm MGAS with the original DIRECT
algorithm proposed in \cite{jones} and its   locally biased
modification LBDirect introduced in \cite{gablo2, gablo3}. These
methods have been chosen for   comparison because they, just as the
MGAS method, do not require the knowledge of the objective function
gradient and work with several Lipschitz constants simultaneously.
The Fortran implementation of the two methods described in
\cite{gablo2} and downloadable from \cite{gablo1} have been used in
both methods.

 To execute
numerical experiments with the algorithm MGAS, we should define its
parameter $\xi$ from (\ref{cond}). In DIRECT (see \cite{jones}),
where a similar parameter is used, the value $\xi$ is related to the
current minimal function value $f_{min}(k)$ and is fixed as follows:
\begin{equation} \label{csi}
\xi = \epsilon |f_{min}(k)|, \hspace{1cm} \epsilon \geq 0.
\end{equation}
The choice of $\epsilon$ between $10^{-3}$ and $10^{-7}$ has
demonstrated good results for DIRECT on a set of test functions (see
\cite{jones}). Since the value $\epsilon = 10^{-4}$ has produced the
most robust results for DIRECT (see, e.g., \cite{gablo2, gablo3,
jones}), exactly this value was used in (\ref{csi}) for DIRECT in
our experiments. The same formula (\ref{csi}) and the same value
$\epsilon = 10^{-4}$ were used in the new algorithm, too.


The series of experiments involves a total of 800 test functions in
the dimensions $N=2,3,4,5$ generated by the GKLS-generator described
in \cite{gaviano} and downloadable  from
http://wwwinfo.deis.unical.it/$\sim$yaro/GKLS.html. More precisely,
eight classes of 100 functions have been considered.
 The generator allows one to construct classes of randomly
 generated
multidimensional and multiextremal test functions with
\textit{known} values of local and global minima and their
locations. Each test class
contains 100 functions and only the following five parameters should be defined by the user: \\
\hspace*{0.5cm} $N$ -- problem dimension; \\
\hspace*{0.5cm} $m$ -- number of local minima; \\
\hspace*{0.5cm} $f^*$ -- value of the global minima; \\
\hspace*{0.5cm} $r^*$ -- radius of the attraction region of the
global minimizer; \\
\hspace*{0.5cm} $d$ -- distance from the global minimizer to the
vertex of the paraboloid.

\begin{figure}[!tbp]
\centering
\includegraphics[scale=.70]{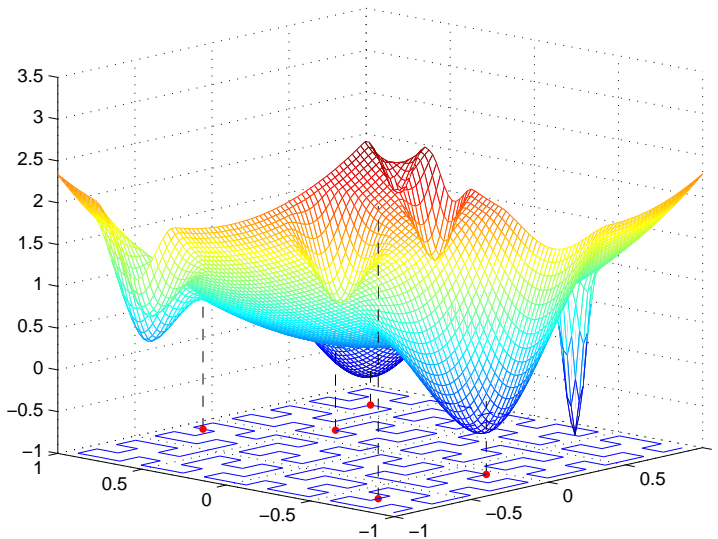}
\hspace{0.0cm} \caption{\emph{A function produced by the
GKLS-generator shown together with a piecewise-linear approximation
to Peano curve used for optimization.}} \label{fig.6}
\end{figure}

 The generator works by constructing in $\Re^N$ a convex
quadratic function, i.e., a paraboloid, systematically distorted
by polynomials. In our numerical experiments we have considered
classes of continuously differentiable test functions with $m=10$
local minima. The global minimum value $f^*$ has been fixed equal
to $-1.0$ for all classes. An example of a function generated by
the GKLS can be seen in Fig.~\ref{fig.6}.

By changing the user-defined parameters, classes with different
properties can be created. For example, a more difficult test class
can be obtained either by decreasing the radius $r^*$ of the
attraction region of the global minimizer or by increasing the
distance, $d$, from the global minimizer to the paraboloid vertex.
In this paper, for each dimension $N=2,3,4,5$, two test classes
where considered: a simple one and a difficult one, see Table
\ref{table1} that describes the classes used in the experiments.
Since the GKLS-generator provides functions with known locations of
global minima, the experiments have been carried out by using the
following stopping criteria.

\begin{table}[t]
\centering \footnotesize
\begin{tabular}{|c|c|c|c|c|c|c|}
\hline
Class  &  Difficulty & N & m & $f^*$ & $d$ & $r^*$  \\
\hline \hline
  1 & Simple & 2 & 10 & -1.0 & 0.90 & 0.20  \\
  2 & Hard   & 2 & 10 & -1.0 & 0.90 & 0.10  \\ \hline
  3 & Simple & 3 & 10 & -1.0 & 0.66 & 0.20  \\
  4 & Hard   & 3 & 10 & -1.0 & 0.90 & 0.20  \\ \hline
  5 & Simple & 4 & 10 & -1.0 & 0.66 & 0.20  \\
  6 & Hard   & 4 & 10 & -1.0 & 0.90 & 0.20  \\ \hline
  7 & Simple & 5 & 10 & -1.0 & 0.90 & 0.40  \\
  8 & Hard   & 5 & 10 & -1.0 & 0.90 & 0.30  \\ \hline
\end{tabular}
\caption {\label{table1} \em Description of 8 classes of test
functions used in  experiments}
\end{table}

 {\bf Stopping criteria.} If $y^*_i$ denotes the global minimizer of the $i$-th function
of the test class, $1\leq i \leq 100$, then the search terminates
either when the maximal number of trials $T_{max}$, equal to
$1000000$, was reached or when a trial point  falls in a ball
$B_i$ having a radius $\rho$ and the center at the global
minimizer of the  $i$-th function of the class, i.e.,
\begin{equation} \label{ball}
 B_i = \{y \in R^N : \|y-y_i^*\| \leq \rho \}, \hspace{5mm} 1\leq i \leq
 100.
\end{equation}
All the methods under comparison can execute $p(k)>1$ trials in the
course of each $k$-th iteration, therefore, when condition
(\ref{ball}) is satisfied at an iteration $k^*$ the number of trials
executed to solve the problem is calculated as   $\sum _{k=1}^{k^*}
p(k)$.   The radius $\rho$ from (\ref{ball}) in the stopping rule
was fixed equal to $0.01\sqrt N $ for classes 1--6 and $0.02\sqrt N
$ for classes 7 and 8.

\begin{table}[t]
\renewcommand\arraystretch{1.2}
\begin{center} \scriptsize
\begin{tabular}{|c|r|r|c|r|r|}\hline
 \multicolumn{3}{|c|}{Class 1}  &\multicolumn{3}{|c|}{Class 5} \\
\cline{1-3} \cline{4-6}  $\eta$ &Average & Maximal & $\eta$ &Average
&Maximal \\ \hline
{$ \bf 10^{-4}$} &{\bf 174.24} & {\bf 565} &$10^{-8}$ & 12174.20 & 171561 \\
$10^{-6}$ &227.60 & 889 & {$ \bf 10^{-10}$} &{\bf 10674.30} &{\bf 95467} \\
$10^{-8}$ &268.98&1279 &$10^{-12}$ &15145.12 &143075 \\
\hline
\end{tabular}
\end{center}
\caption{\em Sensitivity analysis for the parameter $\eta$ from
(\ref{cond1}) } \label{table4}
\end{table}

In order to show the influence  of the parameter $\eta$ introduced
in (\ref{cond1}) on the search,   a sensitivity analysis has been
performed.  Two different classes of test functions have been
considered: class 1 in  $N=2$ and class 5 in   $N=4$ (see
Table~\ref{table1}). Three different values of the parameter $\eta$
were used for each class. In Table~\ref{table4},    the average and
the maximal  number of function evaluations calculated in order to
satisfy the stopping rule for all 100 functions of each class are
reported. Notice that the best results (shown in bold) were obtained
for $\eta = 10^{-4}$ for the class 1, and for $\eta = 10^{-10}$ for
the class 5. In the case of dimension $N=2$, values of $\eta$
smaller than $10^{-4}$ produce an intensification of the search in
subintervals already well-explored. Conversely, when the dimension
$N$ increases the reduced function in one dimension becomes more
oscillating and it is necessary to reduce the value of the parameter
$\eta$. In general, if   a  too small value of the parameter is
applied, the algorithm continues the search in parts of the domain
that were already well-explored during the previous iterations.
Obviously, if
  a too large value of $\eta$ is used it happens that from a certain
iteration onward, no interval is selected  for the subdivision and
the global solution can be lost.

Taking into account results of  the sensitivity analysis, the
following values have been chosen: $\eta=10^{-4}$ for classes 1 and
2, $\eta=10^{-7}$ for the class 3, $\eta=10^{-8}$ for the class 4
and $\eta=10^{-10}$ for the classes 5--8.

In the algorithm MGAS, an $M$-approximation of the Peano curve has
been considered. In particular the level $M$ of the curve must be
chosen taking in mind the constraint $NM<G $, where $N$ is the
dimension of the problem and $G $ is the number of digits in the
mantissa depending on the computer that is used for the
implementation (see [34] for more details). In our experiments we
had $G =52$, thus the value $M=10$ has been used for all the classes
of test functions.

Results of numerical experiments with the eight GKLS test classes
from Table~\ref{table1} are shown in Table~\ref{table2}. The columns
``Average number of trials''  in Table~\ref{table2} report the
average number of trials performed during minimization of the 100
functions from each GKLS class. The simbol ``$>$'' reflects the
situations when not all functions of a class were successfully
minimized by the method under consideration: that is the method
stopped when 1000000 trials had been executed during minimizations
of several functions of this particular test class. In these cases,
the value 1000000 was used in calculations of the average value,
providing in such a way a lower estimate of the average.  The
columns ``Maximal number of trials'' report the maximal number of
trials required for satisfying the stopping rule  for all 100
functions of the class. The notation ``1000000 ($j$)'' means that
after 1000000 function evaluations the method under consideration
was not able to solve $j$ problems.

\begin{table}[t]
\centering \footnotesize
\begin{tabular}{|c|c c c|c c c|}
\hline Class  &  \multicolumn{3}{c|}{Average number of trials}    & \multicolumn{3}{c|}{Maximal number of trials }    \\
\hline &  DIRECT & LBDirect & MGAS & DIRECT & LBDirect & MGAS\\
\hline\hline
  1  & 208.54 &304.28&174.24& 1159&2665&565\\ \hline
  2  &1081.42&1291.70&622.60& 3201&4245&1749 \\ \hline
  3  &1140.68&1893.02&1153.64 & 13369&20779&5267\\ \hline
  4  &$>$42334.36&5245.72&2077.60& 1000000(4)&32603&9809 \\ \hline
  5  &$>$47768.28&21932.94&9961.70 & 1000000(4)&179383&95467\\ \hline
  6  &$>$95908.99&74193.53&21687.76 & 1000000(7)&372633&319493\\ \hline
  7  &$>$33878.09&31955.06&7306.04 & 1000000(3)&146623&36819\\ \hline
  8  &$>$149578.61&$>$93876.77&23460.00  & 1000000(13)&1000000(1)&96287 \\ \hline
\end{tabular}
\caption {\label{table2} \em Results of experiments}
\end{table}

Table \ref{table3} reports the ratio between the maximal (and
average) number of trials performed by DIRECT and LBDirect with
respect to the corresponding number of trials performed by the new
algorithm MGAS. It can be seen from Table \ref{table3} that the
algorithm MGAS outperforms both competitors significantly on the
given test classes.

\begin{table}[h]
\centering \footnotesize
\begin{tabular}{|c|c c|c c|}
\hline Class  &  \multicolumn{2}{c|}{Average number of trials}    & \multicolumn{2}{c|}{Maximal number of trials }    \\
\hline &  DIRECT/MGAS & LBDirect/MGAS & DIRECT/MGAS & LBDirect/MGAS\\
\hline\hline
  1 &1.20& 1.75 & 2.05&4.72\\ \hline
  2  &1.74&2.07 & 1.83&2.43\\ \hline
  3  &0.99&1.64 & 2.54&3.95\\ \hline
  4  &$>$20.38&2.52 & $>$101.95&3.32\\ \hline
  5  &$>$4.80&2.20 & $>$10.47&1.88\\ \hline
  6  &$>$4.42&3.42 & $>$3.13&1.17\\ \hline
  7  &$>$4.64&4.37 & $>$27.16&6.82\\ \hline
  8  &$>$6.38&$>$4.00 & $>$10.39&$>$10.39\\ \hline
\end{tabular}
\caption {\label{table3} \em Speed up obtained by MGAS  with respect
to its competitors}
\end{table}

\begin{figure}
\centering
\begin{minipage}{.450\linewidth}
\centerline{\resizebox{\textwidth}{!}{\includegraphics{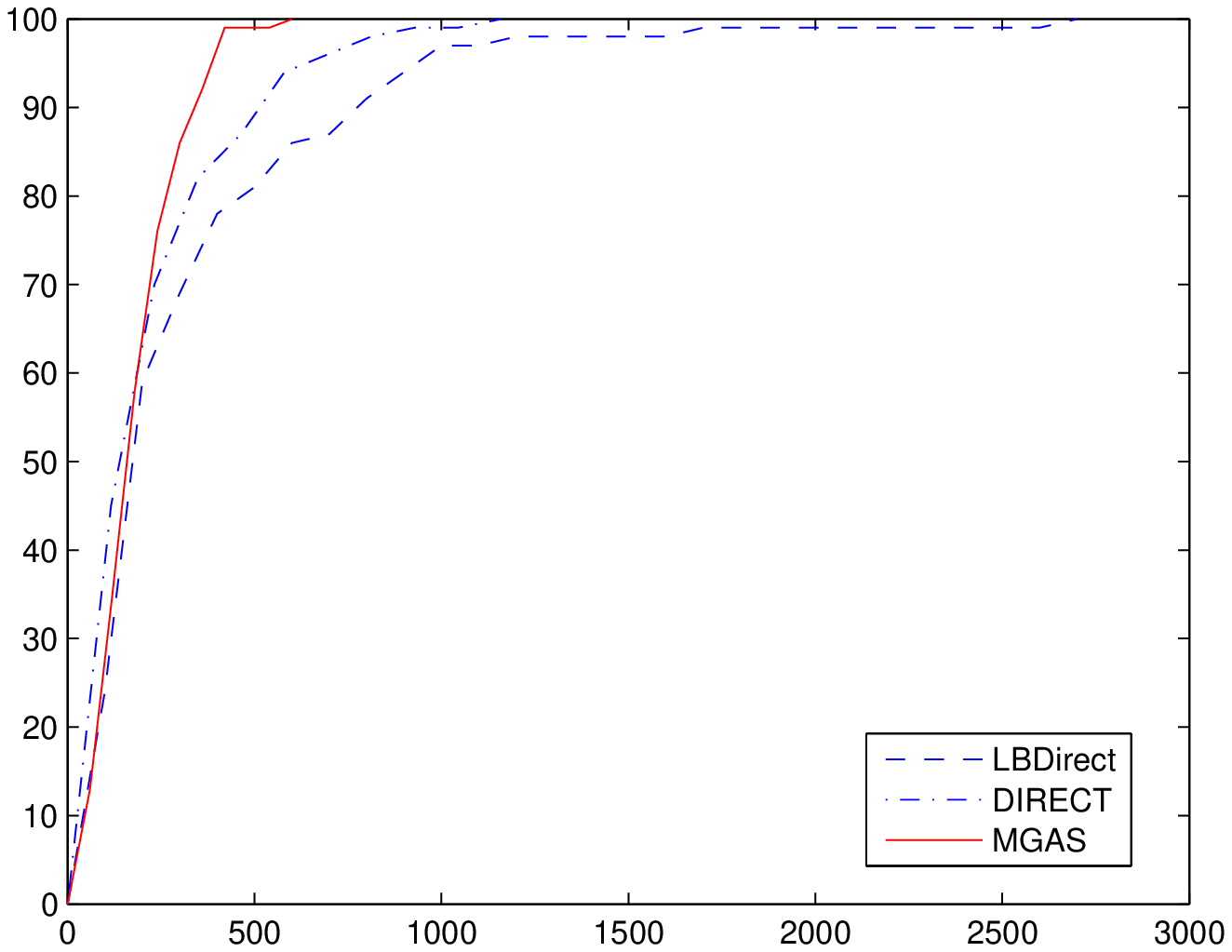}}}
\end{minipage} \hspace{0.1cm}
\begin{minipage}{.450\linewidth}
\centerline{\resizebox{\textwidth}{!}{\includegraphics{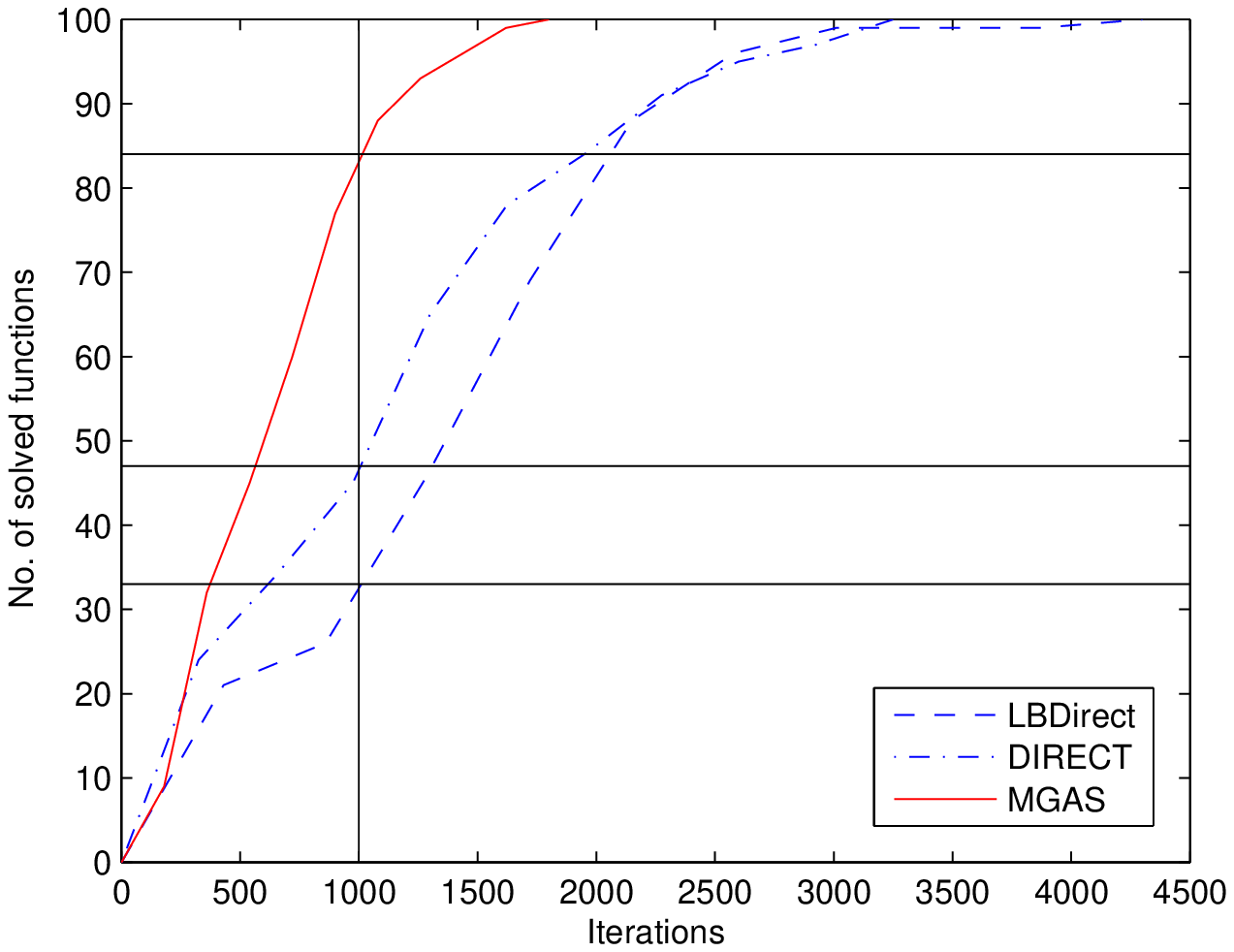}}}
\end{minipage}
\begin{minipage}{.450\linewidth}
\centerline{\resizebox{\textwidth}{!}{\includegraphics{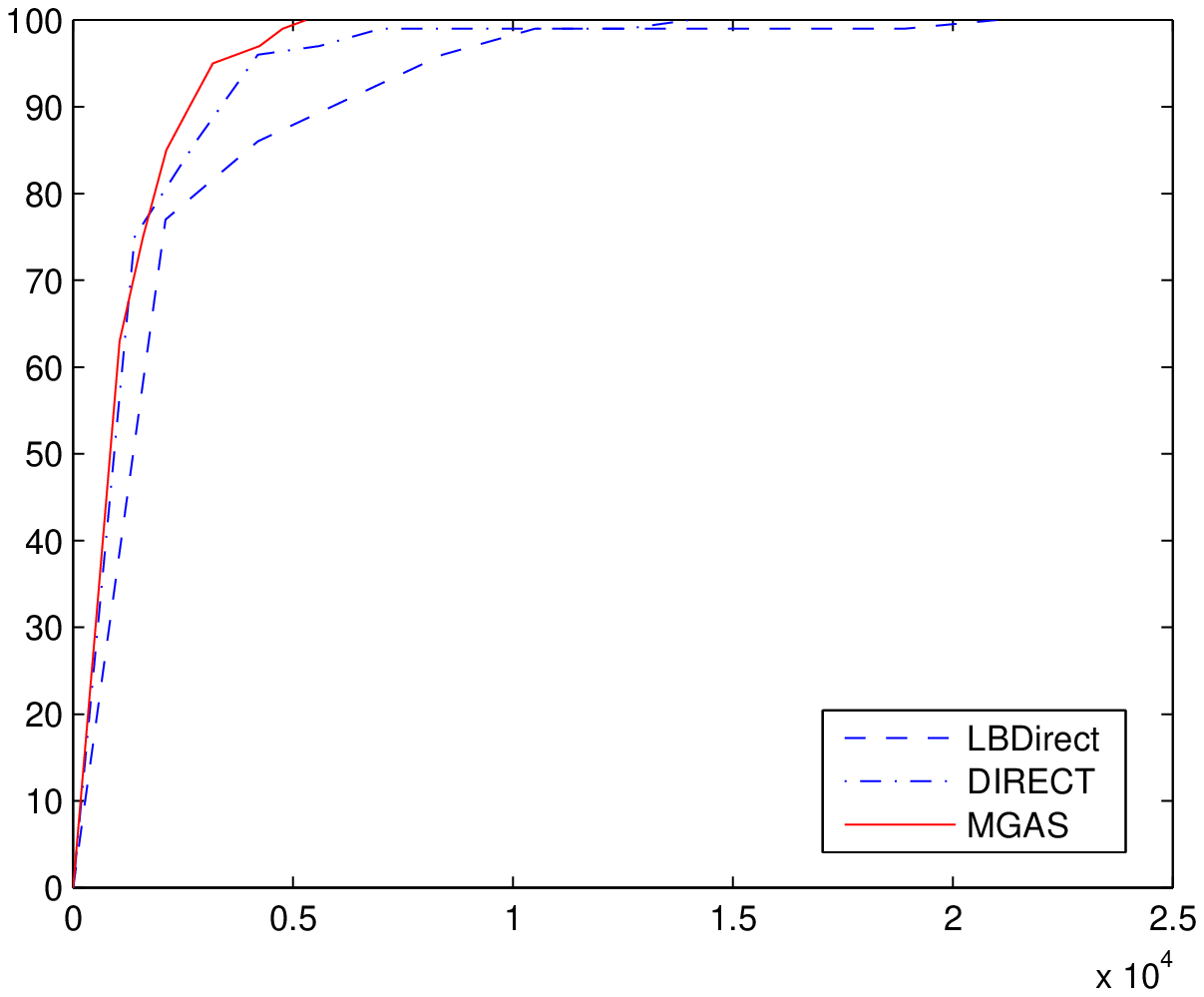}}}
\end{minipage}
\begin{minipage}{.450\linewidth}
\centerline{\resizebox{\textwidth}{!}{\includegraphics{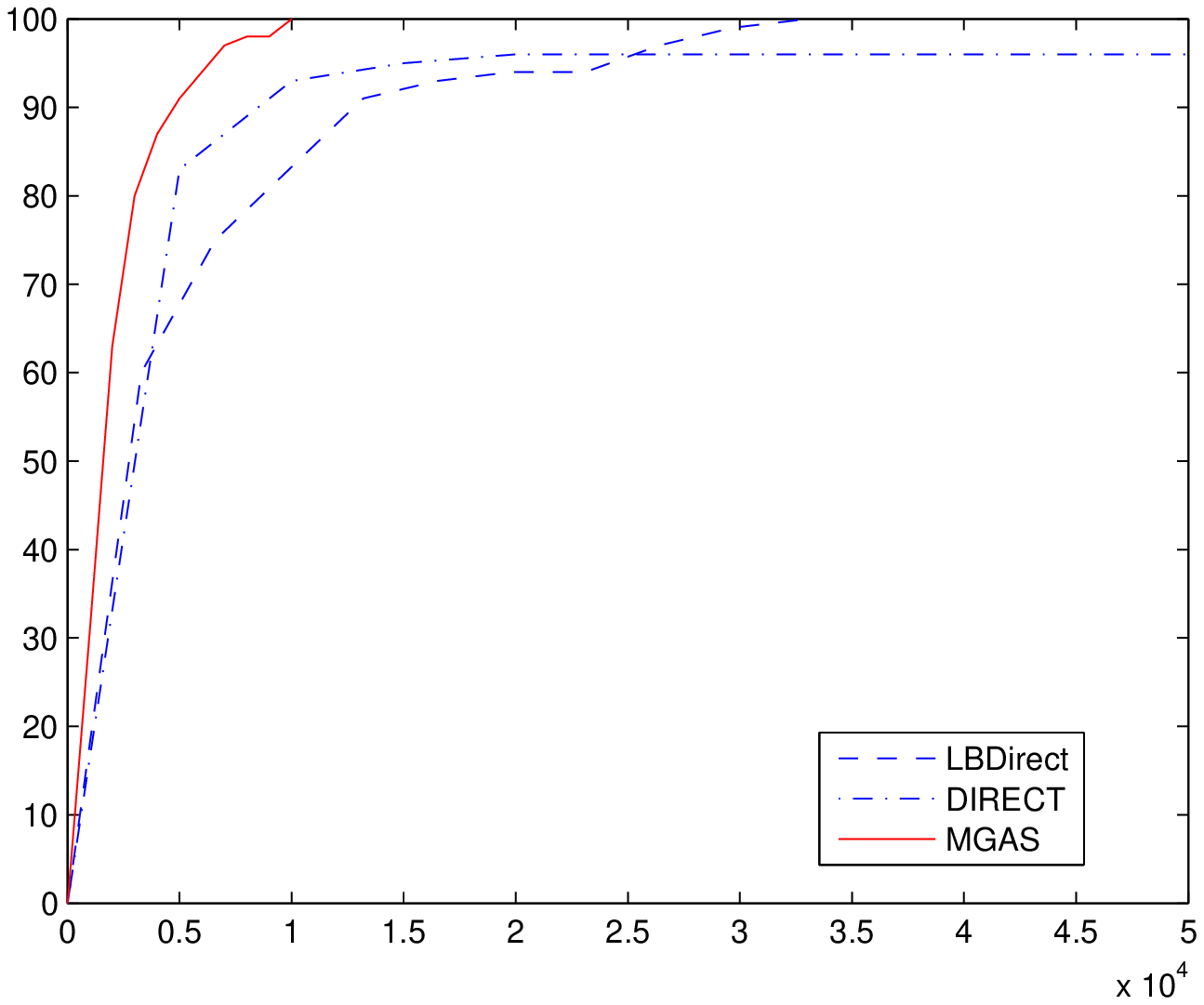}}}
\end{minipage}
\begin{minipage}{.450\linewidth}
\centerline{\resizebox{\textwidth}{!}{\includegraphics{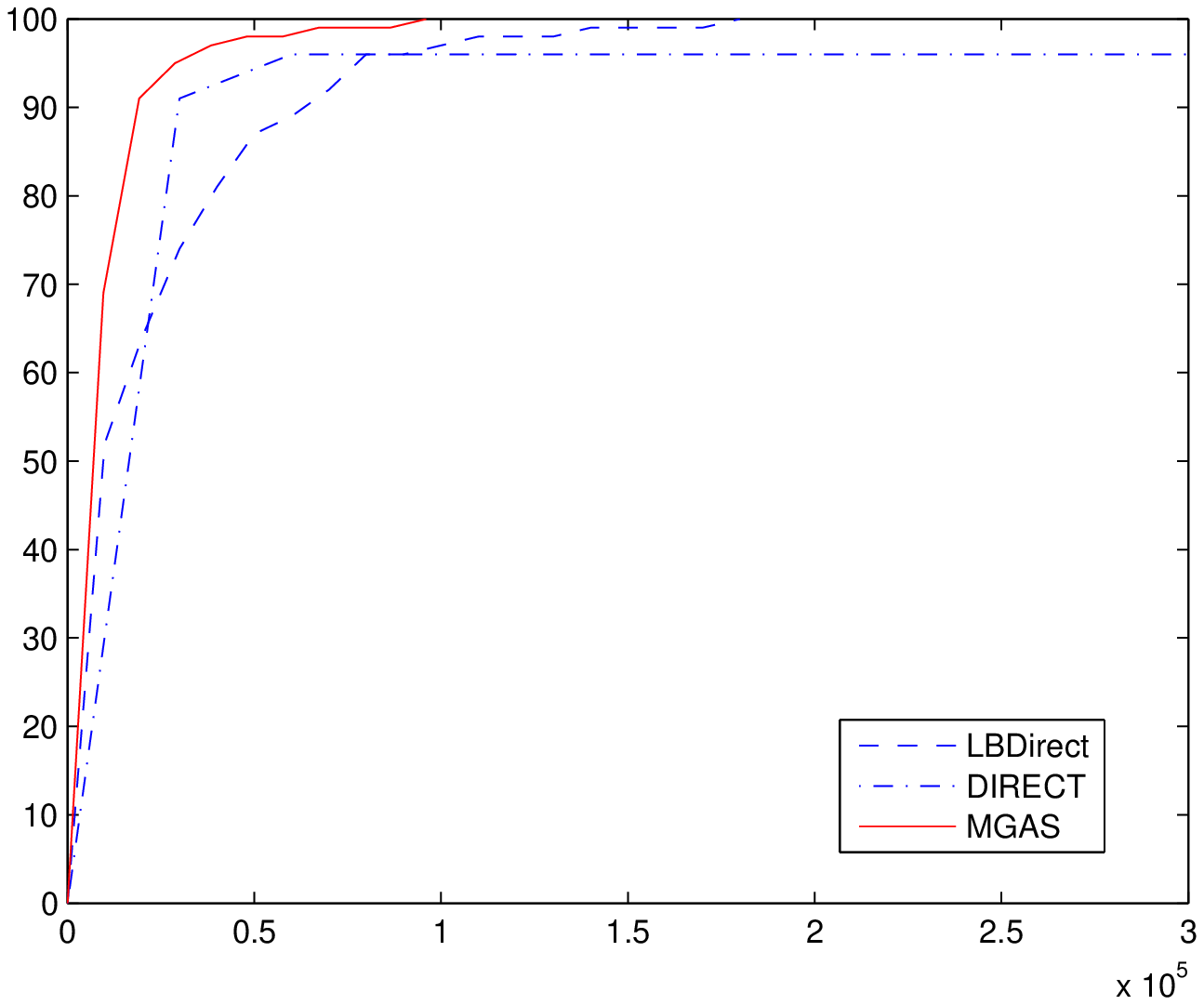}}}
\end{minipage}
\begin{minipage}{.450\linewidth}
\centerline{\resizebox{\textwidth}{!}{\includegraphics{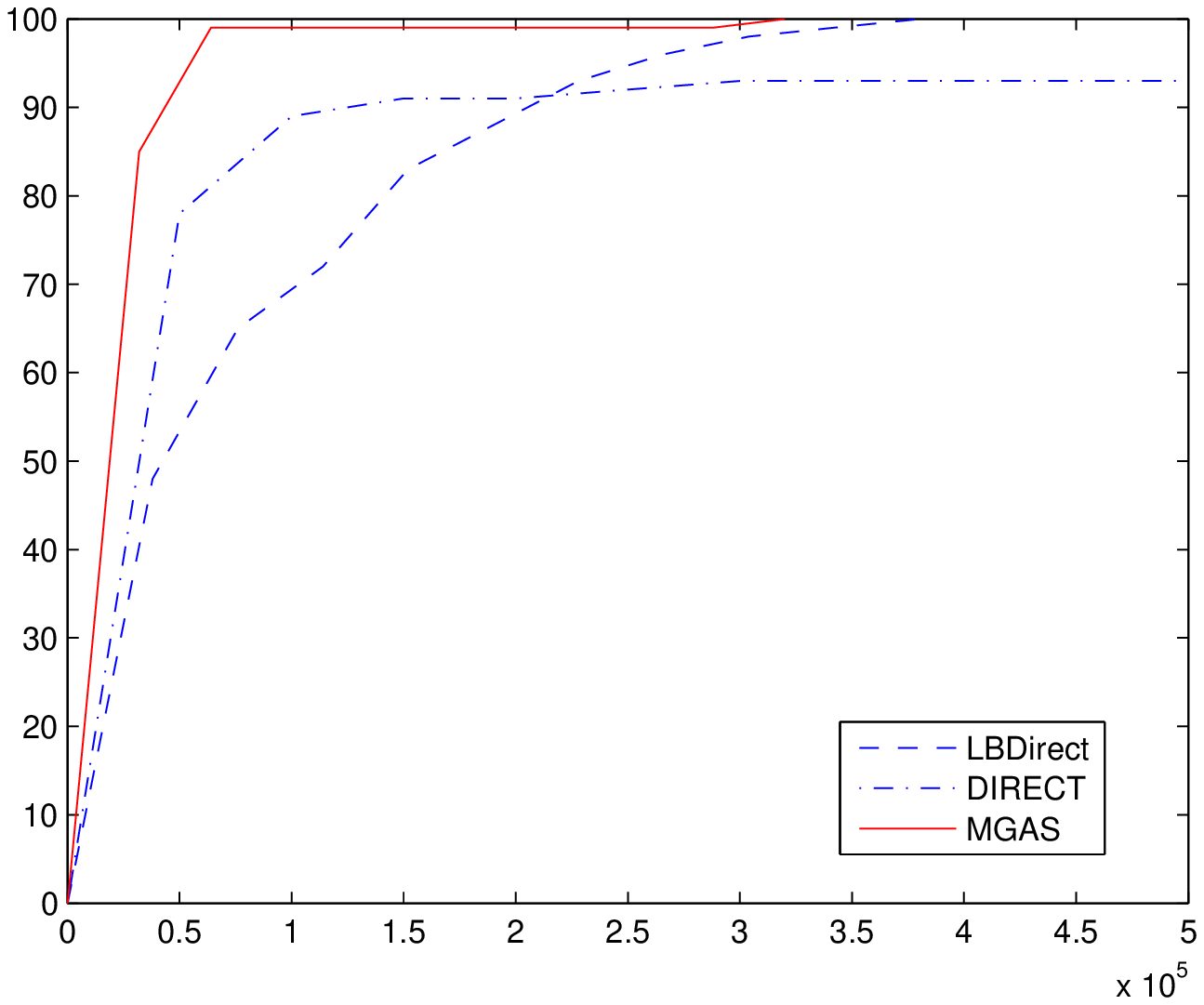}}}
\end{minipage}
\begin{minipage}{.450\linewidth}
\centerline{\resizebox{\textwidth}{!}{\includegraphics{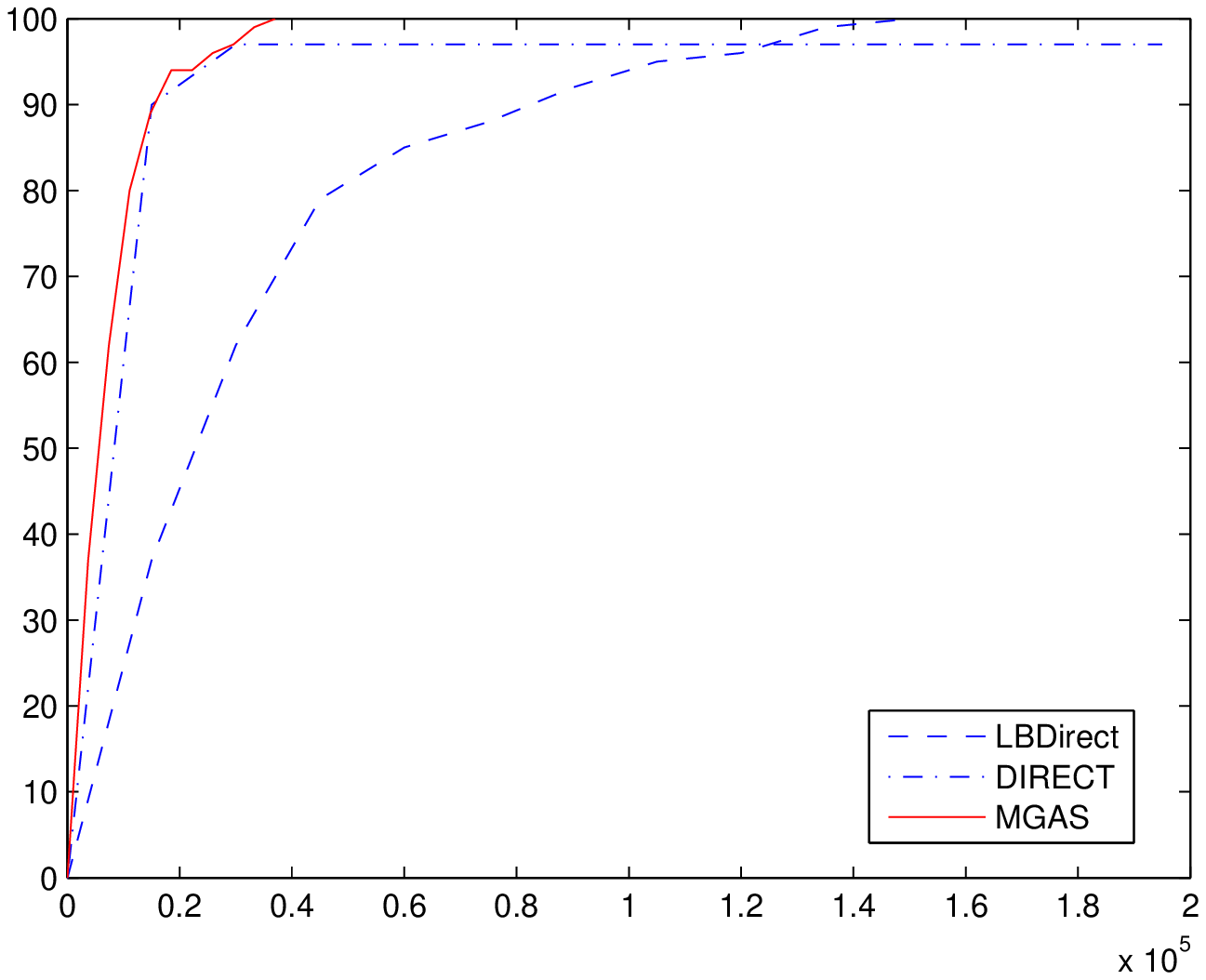}}}
\hspace*{2.5cm}\small{\em classes 1,3,5,7}
\end{minipage}
\begin{minipage}{.450\linewidth}
\centerline{\resizebox{\textwidth}{!}{\includegraphics{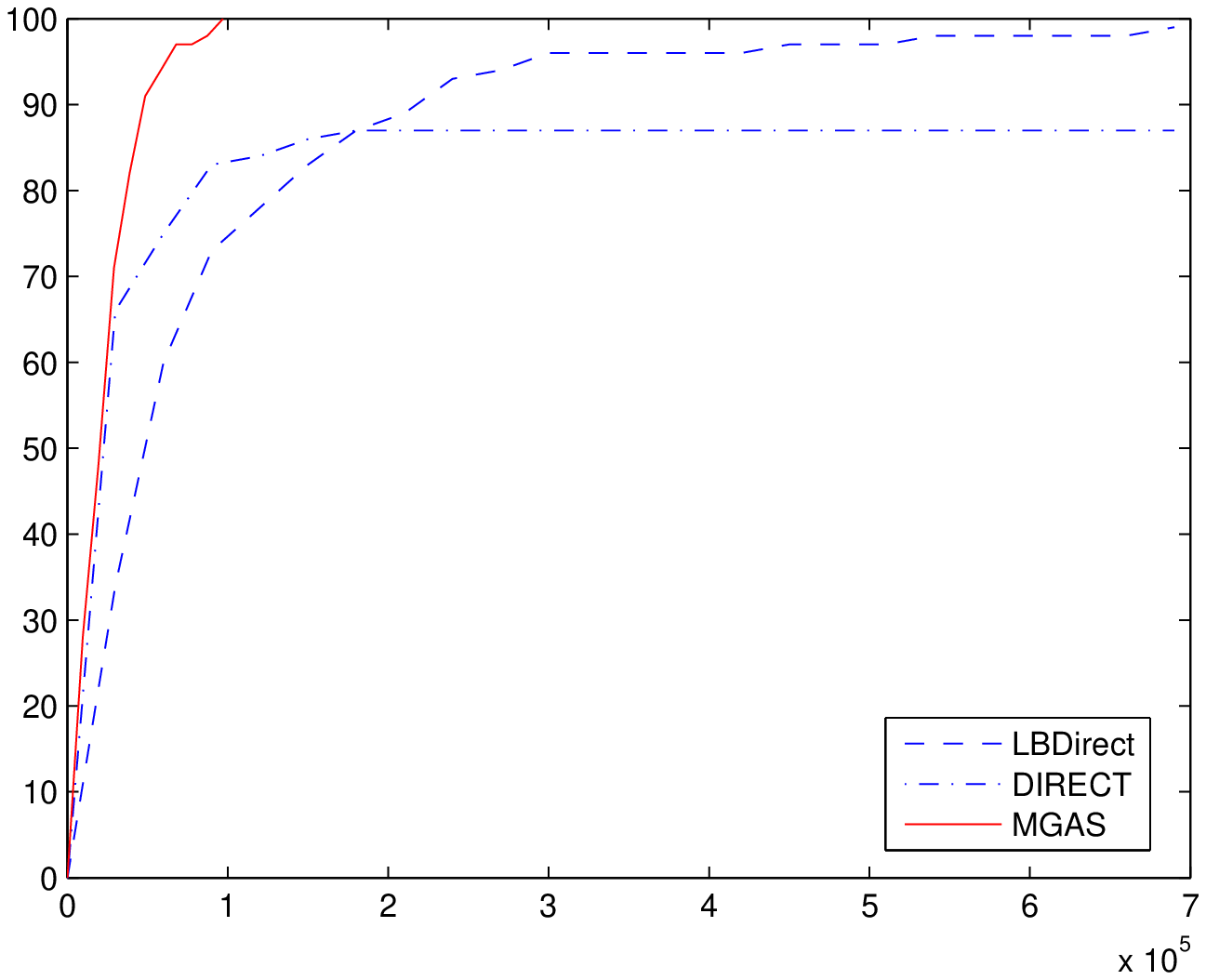}}}
\hspace*{2.5cm}\small{\em classes 2,4,6,8}
\end{minipage}
\caption{\em Operating characteristics for the methods MGAS, DIRECT,
and LBDirect for the eight classes from Table \ref{table1}. The
left-hand column presents results for the simple classes 1,3,5,7,
from top to bottom. The right-hand  column shows results for the
difficult classes 2,4,6,8.} \label{fig.7}
\end{figure}

Figure~\ref{fig.7} shows a comparison of the three methods  using
the so called \textit{operating characteristics} introduced in 1978
in~\cite{I_Grishagin1978} (see, e.g.,~\cite{yaro} for their
English-language description). These characteristics show very well
the performance of   algorithms under the comparison for each class
of test functions. On the horizontal axis  we have the number of
function evaluations and the vertical coordinate of each curve shows
how many problems have been solved by one or another method after
executing the number of function evaluations corresponding to the
horizontal coordinate. For instance, the first graph in the
right-hand column ($N=2$, class 2) shows that after 1000 function
evaluations the LBDirect has found the global solution at 33
problems, DIRECT at 47 problems, and the MGAS at 84 problems. Thus,
the behavior of an algorithm is better if its characteristic is
higher than characteristics of its competitors. In
Figure~\ref{fig.7}, the left-hand column of characteristics, the
behavior of algorithms MGAS, DIRECT, and LBDirect on the classes 1,
3, 5, and 7 is shown. The right-hand column presents the situation
when  the more difficult classes 2, 4, 6, and 8 have been used.

\section{A brief conclusion}

The global optimization problem  of a multi-dimensional,
non-differentiable, and multiextremal function has been considered
in this paper. It was supposed that the objective function   can be
given as a `black-box' and the only available information is that it
satisfies the Lipschitz condition with an unknown Lipschitz constant
over the search region being a hyperinterval in $\Re^N$.

A new deterministic global optimization algorithm called  MGAS has
been proposed. It uses   the following two ideas: the MGAS applies
numerical approximations to space-filling curves to reduce the
original Lipschitz multi-dimensional problem to a univariate one
satisfying the H\"{o}lder condition; the MGAS  at each iteration
uses a new geometric technique working with a number of possible
H\"{o}lder constants chosen from a set of values varying from zero
to infinity evolving so ideas of the popular DIRECT method to the
field of H\"{o}lder global optimization. Convergence conditions of
the MGAS have been established. Numerical experiments carried out on
800 of test functions generated randomly have been executed.

It can be seen from the numerical experiments that the new
algorithm shows quite a promising performance in comparison with
its competitors. Moreover, the advantage of the new technique
becomes more pronounced for harder problems.

\end{document}